\documentclass[10pt]{article}
\usepackage{amsmath}
\usepackage{latexsym}
\usepackage{amssymb}
\usepackage{amsfonts}
\usepackage{amsthm}
\title{HOMOLOGICAL SOLUTION OF THE \\ RIEMANN/LANCZOS AND WEYL/LANCZOS \\
PROBLEMS IN ARBITRARY DIMENSION }
\author{J.-F. POMMARET \\ CERMICS, Ecole des Ponts ParisTech, France  \\
( jean-francois.pommaret@wanadoo.fr)  }
\date{  }
\begin{document}
\maketitle

\noindent
{\bf ABSTRACT}   \\

When ${\cal{D}}$ {is a linear ordinary or partial differential operator  of any order, a {\it direct problem} is to look for an operator ${\cal{D}}_1$ generating the {\it compatibility conditions} (CC) ${\cal{D}}_1\eta=0$ of ${\cal{D}}\xi=\eta$. Conversely, when ${\cal{D}}_1$ is given, an {\it inverse problem} is to look for an operator $\cal{D}$  such that its CC are generated by ${\cal{D}}_1$ and we shall say that ${\cal{D}}_1$ is {\it parametrized} by $\cal{D}$. We may thus construct a {\it formally exact} differential sequence with successive operators ${\cal{D}},{\cal{D}}_1,{\cal{D}}_2, ...$, where each operator is parametrizing the next one. However, this sequence may not be {\it strictly exact} in the sense that certain operators may be neither involutive nor even formally integrable. Introducing the {\it formal adjoint} $ad( \,\,)$, we have ${\cal{D}}_i\circ {\cal{D}}_{i-1}=0 \Rightarrow ad({\cal{D}}_{i-1}) \circ ad({\cal{D}}_i)=0$ but $ad({\cal{D}}_{i-1})$ may not generate {\it all} the CC of $ad({\cal{D}}_i)$.   When $D=K[d_1,...,d_n]=K[d]$ is the (non-commutative) ring of differential operators with coefficients in a differential field $K$, it gives rise by residue to a {\it differential module} $M$ over $D$. The homological {\it extension modules} $ext^i(M)=ext^i_D(M,D)$ with $ext^0(M)=hom_D(M,D)$  only depend on $M$ and are measuring the above gaps, independently of the previous differential sequence.   \\
The purpose of this concise but technical paper is to compute them for certain Lie operators involved in the formal theory of Lie pseudogroups in arbitrary dimension $n$. In particular, we prove for the first time that the extension modules highly depend on the Vessiot {\it structure constants} $c$. When one is dealing with a Lie group of transformations or, equivalently, when ${\cal{D}}$ is a Lie operator of finite type, then we shall prove that $ext^i(M)=0, \forall 0\leq i \leq n-1$. It will follow that the {\it Riemann-Lanczos} and {\it Weyl-Lanczos} problems just amount to prove such a result for $i=2$ and arbitrary $n$ when ${\cal{D}}$ is the {\it Killing} or {\it conformal Killing} operator. We finally prove that ${ext}^i(M)=0, \forall i\geq 1$ for the Lie operator of infinitesimal contact transformations with arbitrary $n=2p+1$. Most of these new results have been checked by means of computer algebra. \\

\vspace{1cm}


\newpage

\noindent
{\bf 1) INTRODUCTION}  \\

\noindent
{\it Ordinary differential (OD) control theory} studies input/output relations
defined by systems of ordinary differential (OD) equations. In this case, with independant variable $x=t=time$, dependant variables 
$y=(y^1, ... ,y^m)$ and standard notation $d=d_x=d_t$, if a control system is defined by input/output ($u\leftrightarrow y$) Kalman relations 
$d_xy=Ay+Bu $, this system is ``{\it controllable}'' if  and only if $rk(B,AB,...,A^{m-1}B)=m$  ([11],[29]). However, and despite many attempts, such a definition seems purely artificial. Now, let us consider the system of two OD equations  for three unknowns where $a(x)$ is a variable parameter:  \\
\[  d_xy^1-a(x)y^2- d_xy^3=0 \hspace{1cm}, \hspace{1cm}  y^1- d_xy^2+d_xy^3=0,  \]
We let the reader check that, whether we choose $y^2$ or $y^3$ as input, we get two quite different systems
in Kalman form, though both are controllable if and only if $a\neq 0$ and $a\neq 1$ whenever $a$ is constant but nothing can be said when $a=a(x)$ is no longer a constant.\\
{\it Two problems are raised at once}.\\
First of all, if the derivatives of the inputs do appear in the control system, for example in the SISO system $\dot{x}-\dot{u}=0$, not a word is
left from the original definition of controllability which is only valid for systems in ``{\it Kalman form}''.
Secondly, we understand from the above example that {\it controllability must be a structural property of a control system}, 
neither depending on the choice of the inputs and outputs among the system variables, nor even on the presentation of the system (change of the
variables eventually leading to change the order of the system). For example, using the second ODE to get $y^1$ ans substituting into the first, 
we should get the second order ODE:  \\
\[                d_{xx}y^2 - a(x)y^2 - d_{xx}y^3 -  d_xy^3=0   \]
More generally, ``{\it partial differential (PD) control theory}'' will study input / output relations defined by systems of partial differential (PD)
equations.  At first sight it does not seem that we have any way to generalize the Kalman form and not a word of the preceding
approach is left as, in most cases, the number of arbitrary parametric derivatives could be infinite. We also
understand that a good definition of controllability should also be valid for control systems with variable coefficients.  \\

Keeping aside these problems for the moment, let us now turn to the formal theory of systems of OD or PD equations.\\

In 1920, M. Janet provided an effective algorithm for looking at the formal (power series) solutions of systems of ordinary differential (OD) or 
partial differential (PD) equations ([10]). The interesting point, in the approach of Janet, is that it also allows to determine 
the {\it compatibility conditions} ${\cal D}_1\eta=0$ for solving {\it formally} inhomogeneous systems of the form $\cal D \xi=\eta$ when 
$\cal D$ is an OD or PD operator and $\xi, \eta$ certain functions. Similarly, one can also determine the compatibility conditions ${\cal D}_2\zeta=0$ 
for solving ${\cal D}_1 \eta=\zeta$, and so on. With no loss of generality, this construction of a ``{\it differential sequence}'' can be done in such a canonical way that we successively obtain ${\cal D}_1, {\cal D}_2, ..., {\cal D}_n$ from $\cal D$ and ${\cal D}_n$ is surjective when $n$ is the number of independent variables. It must be noticed that this important result has been just provided as a footnote in ([10]).   \\

With no reference to the above work that he was ignoring, D.C. Spencer developed, from 1965 to 1975, the {\it formal theory} of 
systems of PD equations by relating the preceding results to {\it homological algebra} and {\it jet theory} ([23],[27],[44]). However, this tool has almost been ignored by mathematicians and, ``{\it a fortiori}'', by engineers or even physicists. Therefore, it becomes clear that the module 
theoretic counterpart, today known as ``{\it algebraic analysis}'', which has been pioneered around 1970 by V.P. Palamodov for the constant 
coefficient case ([22]) and by M. Kashiwara ([12]) for the variable coefficient case, as it heavily depends on the previous 
difficult work and looks like even more abstract, has been {\it totally ignored} within the range of any application before 1990, 
when U. Oberst revealed its importance for control theory ([19], Compare to [28]).    \\

As we always use to say, the difficulty in studying differential ideals or differential modules is not of an algebraic nature but rather of a
differential geometric nature. This is the reason for which the study of algebraic analysis is at once touching delicate points of jet theory, the
main ones being {\it formal integrability} and {\it formal duality} with their strange relationship. We shall explain these concepts on a few tricky motivating examples and invite the reader to try to discover {\it now} by himself that, when $a=a(x)$ in the last example, the " {\it built-in} " controllability condition becomes $d_xa+ a^2 - a\neq 0$, in order to understand the powerfulness of these new methods. \\

\noindent
{\bf DEFINITION 1.1}: ({\it Formal integrability}) We say that a system of OD or PD equations of order $q$ is {\it formally integrable} if, whenever we differentiate all the equations $r+1$ times, then it does not bring more equations of order $q+r$ than if we were differentiating these equations only $r$ times, $\forall r\geq 0$. Similarly, an operator will be said to be formally integrable if the corresponding system is formally integrable. An operator/system will be said to be {\it finite type} if all the derivatives of unknowns are known at a certain order from the ones of lower order. An operator/system will be said to be {\it involutive} if it is formally integrable and the Spencer $\delta$-cohomology of its symbol vanish.   \\

As this first delicate point is not very familiar to the computer algebra community, we shall carefully distinguish various types of differential sequences (See [23],[26],[33] for the details):  \\
\noindent
$\bullet$ A differential sequence is said to be {\it formally exact} if each operator generates the CC of the previous one.  \\
\noindent
$\bullet$ A formally exact differential sequence is said to be {\it strictly exact} if each operator is formally integrable.  \\
\noindent
$\bullet$ A strictly exact differential sequence is said to be {\it involutive} if each operator is involutive.  \\

We have therefore the classification of differential sequences: \\

\hspace*{2cm}  INVOLUTIVE  $\subset$ STRICTLY EXACT   $\subset $  FORMALLY EXACT         \\

\noindent
{\bf EXAMPLE 1.2}: The Poincar\'{e}, Janet and "second" Spencer sequences are involutive but the "first" Spencer sequence is not even strictly exact because the Spencer operator is not involutive as we have ${\partial}_i({\partial}_jf^k - f^k_j)- {\partial}_j({\partial}_if^k-f^k_i)={\partial}_jf^k_i - {\partial}_if^k_j$.   \\

\noindent
{\bf DEFINITION 1.3}: ({\it Formal adjoint}) If ${\cal{D}}:\xi \rightarrow \eta$ is a given differential operator of order $q$, then the {\it formal adjoint} operator $ad({\cal{D}}): \mu \rightarrow \nu$ is an operator of the same order $q$ defined by multiplying the equations ${\cal{D}}\xi=\eta$ by test functions $\lambda$, summing and integrating by parts along the following formula:   \\
\[       < \mu, {\cal{D}}\xi> = < ad({\cal{D}}) \mu , \xi>  \,\,+ \,\, div( ... )   \]
We have $ad(ad({\cal{D}}))={\cal{D}}$ and most of the results of this paper will come from the fact that an operator may be formally integrable while its formal adjoint is not at all formally integrable.  \\

\noindent
{\bf DEFINITION 1.4}: From now on, all the operators considered will have coefficients in a differential field $K$ with derivations 
$({\partial}_1, ... , {\partial}_n)$, for example in $\mathbb{Q}(a)$ or $\mathbb{Q}(x^1,x^2,x^3)$ and we shall denote by 
$D=K[d_1, ...,d_n]=K[d]$ the (non-commutative) ring of partial differential operators with coefficients in $K$. The following result will be quite useful for applications ([27],[28]):  \\

\noindent
{\bf THEOREM 1.5}: Denoting simply by $rk_D({\cal{D}})$ the {\it differential rank} of ${\cal{D}}$ over $D$, that is the maximum number of equations differentially independent over $D$, we have $rk_D({\cal{D}})=rk_D(ad({\cal{D}}))$.  \\

\noindent
{\bf MOTIVATING EXAMPLE 1.6}: With two independent variables $(x^1,x^2)$, one unknown $y$ and standard notations with $K=\mathbb{Q}(x^1,x^2)$, we consider the differential operator ${\cal{D}}$ defined by the following third order system of PD equations with second member $(u, v)$:
\[    \left\{   \begin{array}{rcll}
Py & \equiv & d_{222}y+x^2y &=u\\
Qy & \equiv & d_2y+d_1y     &=v
\end{array}   \right.   \]
We check the identity $QP-PQ\equiv 1$ and obtain easily the splitting operator:    \\
\[ y=Qu-Pv=d_2u+d_1u-d_{222}v-x^2v  \]
Substituting in the previous PD equations, we should obtain the generating
$6^{th}$-order compatibility conditions for $(u, v)$ in the form: 
\[
\left\{
\begin{array}{rll}
A\equiv &PQu-P^2v-u & =0\\
B\equiv &Q^2u-QPv-v & =0
\end{array}
\right.  \]
These two compatibility conditions are differentially dependent as
we check at once $QA-PB\equiv 0$.\\
Finally, setting $u=0, v=0$, we notice that the preceding homogeneous 
system can be written in the form ${\cal D}y=0$ and admits the only 
solution $y=0$. Of course, the system/operator is not formally integrable because, differentiating the second PD equation with respect to $d_{22}$, then substracting the first equation, we discover that the system is finite type at order  $3$ with:  \\
\[  y_{222} + x^2 y=0, \,\,\,  y_{122} - x^2y=0, \,\,\,   y_{112}-x^2y=0,    \,\,\,     y_{111}+x^2y=0     \]
and even at order $0$ with $y=0$. Using the standard notation $Dy^1 + ...+ Dy^m=Dy\simeq D^m$ for free differential modules, we obtain the locally and formally exact following split {\it differential sequence} which is far from being strictly exact:    \\
\[        0  \rightarrow  D \underset{3}{\stackrel{{\cal{D}}_2}{\longrightarrow}} D^2 \stackrel{{\cal{D}}_1}{\underset{6}{\longrightarrow}} D^2 \underset{3}{\stackrel{{\cal{D}}}{\longrightarrow}}  D   \rightarrow  0    \]
where we have indicated the order of an operator under its arrow. The successive orders are $(3,6,3)$ but we shall see later on examples with $(1,3,1), (1,2,2,1), (1,2,1,2,1)$ and so on, contrary to the situation of an involutive differential sequence where the orders are $(q,1,...,1)$ ([23],[26],[28]).  \\

\noindent
{\bf MOTIVATING EXAMPLE 1.7}: (See [26] and [29] for more details). With three
independent variables $(x^1,x^2,x^3)$ and one unknown $y$, let us
consider the second order system with variable coefficients:
\[
\left\{
\begin{array}{rll}
Py\equiv & d_{33}y-x^2d_{11}y &=u \\
Qy\equiv &d_{22}y             &=v
\end{array}
\right.   \]
Introducing as before:
\[  d_{112}y=\frac{1}{2}(d_{33}v-x^2d_{11}v-d_{22}u)=w \]
we finally get the two following generating compatibility conditions:
\[
\left\{
\begin{array}{rll}
A\equiv &d_{233}v-x^2d_{112}v-3d_{11}v-d_{222}u & =0\\
B\equiv &d_{3333}w-2x^2d_{1133}w+(x^2)^2d_{1111}w-d_{11233}u+
x^2d_{11112}u-d_{1111}u &=0
\end{array}
\right.   \]
These two compatibility conditions of respective orders 3 and 6 are
differentially dependent as one checks at once through computer algebra:
\[ d_{3333}A-2x^2d_{1133}A+(x^2)^2d_{1111}A-2d_2B=0  \]
The space of solutions of the system ${\cal D}y=0$ can be described by polynomials in $(x^1, x^2, x^3)$ and has finite 
dimension 12 over the field of constants $C=\{ a\in K\mid {\partial}_i a=0, \forall i=1,...,n\}\subset K$. We have the formally exact sequence:  \\
\[           0 \rightarrow D \underset{4}{\rightarrow }  D^2 \underset{6}{\rightarrow} D^2 \underset{2}{\rightarrow} D  \rightarrow M  \rightarrow 0   \]

\noindent
{\bf MOTIVATING EXAMPLE 1.8}: Again with two independent variables $(x^1,x^2)$, one unknown $y$ and $K=\mathbb{Q}$, let us consider the following second order system with constant coefficients:
\[\left\{    \begin{array}{rll}
Py\equiv &d_{22}y   &=u\\
Qy\equiv &d_{12}y-y &=v
\end{array}
\right.   \]
We obtain at once the splitting operator:   \\
\[  y=d_{11}u-d_{12}v-v\]
and could hope to obtain as before the generating $4^{th}$-order generating compatibility conditions by substitution, that is to say:
\[ \left\{  \begin{array}{rll}
A\equiv &d_{1122}u-d_{1222}v-d_{22}v-u&=0\\
B\equiv &d_{1112}u-d_{11}u-d_{1122}v &=0
\end{array}   \right.  \]
However, {\it in this particular case}, we notice that $PQ - QP=0$ and there is an unexpected {\it unique second order} generating 
compatibility condition ${\cal{D}}'_1$ of the form:
\[  C\equiv d_{12}u-u-d_{22}v=0 \]
as we have indeed $A\equiv d_{12}C+C$ and $B\equiv d_{11}C$, a result 
leading to $C\equiv d_{22}B-d_{12}A+A$. Accordingly, the systems $A=0, B=0$
on one side and $C=0$ on the other side are completely different though 
they have the same solutions in $u, v$. Accordingly, defining the differential module $M$ by residue with a caninical projection $p$, we have the two completely different resolutions: \\
\[   0 \rightarrow D \stackrel{{\cal{D}}_2}{\underset{2}{\longrightarrow}} D^2 \stackrel{{\cal{D}}_1}{\underset{4}{\longrightarrow}} D^2 \stackrel{{\cal{D}}}{\underset{2}{\longrightarrow}} D  \rightarrow 0 \hspace{5mm} \Leftrightarrow   \hspace{5mm}  
          0 \rightarrow D \stackrel{{\cal{D}}'_1}{\underset{2}{\longrightarrow}}  D^2   \stackrel{{\cal{D}}}{\underset{2}{\longrightarrow}}  D   \rightarrow 0  \]
which are again far from being strictly exact with $M=0$ in both sequences.  \\

After recalling the mathematical tools needed in Section $2$, we study a few specific situations in section $3$, paying a particular attention to 
the Riemann/Lanczos and Weyl/Lanczos problems. The author thanks Prof. Lars Andersson (AEI, Potsdam) for having pointed out to him the interest of using Algebraic Analysis in order to revisit the work of Lanczos in a modern setting.  \\

\noindent
{\bf 2) MATHEMATICAL TOOLS}\\

\noindent
{\bf A) Differential sequences}  \\

\noindent
In view of the many cases illustrated by the preceding similar examples, it becomes clear that there is a need for classifying the
properties of systems of PD equations in a way that does not depend on their presentations and {\it this is the purpose of algebraic analysis} along with 
the scheme:   \\
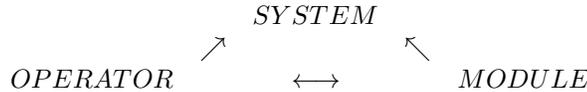
\[\begin{array}{rrrclll}
        &   & SYSTEM &    &  \\    
         &    \nearrow & & \nwarrow   &  \\
   OPERATOR & &  \longleftrightarrow \hspace{5mm} & &  MODULE 
\end{array} \]
in order to show that certain concepts, which are clear in one framework, may become quite obscure in the others and conversely, like the formal integrability and torsion concepts for example.\\

When $E$ is a vector bundle over $X$, having a system of order $q$ on $E$, say $R_q\subset J_q(E)$, we can introduce 
the canonical projection $\Phi:J_q(E)\longrightarrow J_q(E)/R_q=F$ and 
define a linear differential operator 
${\cal{D}}:E\longrightarrow F:\xi(x)\longrightarrow 
{\eta}^{\tau}(x)=a^{\tau\mu}_k(x){\partial}_{\mu}{\xi}^k(x)$. When $\cal{D}$
is given, the compatibility conditions for solving $\cal{D}\xi=\eta$ can be 
described in operator form by ${\cal{D}}_1\eta=0$ and so on. In general 
(see the preceding examples), if a system is not formally integrable, it 
is possible to obtain a formally integrable system, having the 
same solutions, by ``{\it saturating}'' conveniently the given PD equations 
through the adjunction of new PD equations obtained by various prolongations 
and {\it such a procedure must absolutely be done before looking for 
the compatibility conditions}.\\

In order to study differential modules, for simplicity we shall forget
about changes of coordinates and consider trivial bundles. If $K$ is a 
{\it differential field} with $n$ commuting derivations 
${\partial}_1,...,{\partial}_n$ (Say $\mathbb{Q}, \mathbb{Q}(a)$ 
or $\mathbb{Q}(x^1,...,x^n)$ in the previous examples), we denote by $C=cst(K)$ the 
subfield of {\it constants} of $K$, that is the set of elements killed by 
the $n$ derivations (Say $\mathbb{Q}$ or $\mathbb{Q}(a)$ when $a\in C$ in the previous
examples). If $d_1,...,d_n$ are {\it formal derivatives} (pure symbols in
computer algebra packages !) which are only supposed to satisfy 
$d_ia=ad_i+{\partial}_ia$ in the operator sense for any $a\in K$, we may 
consider the (non-commutative) ring $D=K[d_1,...,d_n]$ of differential 
operators with coefficients in $K$. If now $y=(y^1,...,y^m)$ is a set of 
differential indeterminates, we let $D$ act formally on $y$ by setting 
$d_{\mu}y^k=y^k_{\mu}$ and set $Dy=Dy^1+...+Dy^m$. We may also set ${\Phi}^{\tau}\equiv a^{\tau\mu}_ky^k_{\mu}\stackrel{d_i}{\longrightarrow} d_i{\Phi}^{\tau}\equiv a^{\tau\mu}_ky^k_{\mu+1_i}+{\partial}_ia^{\tau\mu}_ky^k_{\mu}$ for $\tau=1,...,p$. Denoting 
simply by $D{\cal{D}}y$ the subdifferential module generated by all the given OD or PD 
equations and all their formal derivatives, we may finally  introduce the $D$-module $M=Dy/D{\cal D}y$ by residue. Here we recall that 
$M$ is a {\it module} over a ring $A$ or an $A$-module if $\forall a\in A, \forall m,n\in M \Rightarrow am, m+n\in M$. We may introduce as usual the {\it torsion} submodule $t(M)=\{ m\in M\mid \exists 0\neq a\in A, am=0 \}$ and we say that $M$ is a torsion module if $t(M)=M$ or that $M$ is torsion-free if $t(M)=0$.\\

\noindent
{\bf EXAMPLE 2.A.1}: In Examples 1.6 and 1.8, we get $M=0$. In Example 1.7, with
$K=\mathbb{Q}(x^1, x^2, x^3)$, we know from ([26], Introduction and Example III.D.1 and [29])
that $M$ is a finite dimensional vector space over $K$ with $dim_K(M)=12$. \\

{\it It is not evident at all} to exhibit the link existing between these two approaches and we proceed as follows. First of all, the ring $D$ is filtred by the order of the operators and we have the filtration or {\it inductive limit} $0=D_{-1}\subset D_0\subset D_1\subset  ... \subset D_q \subset ... \subset D_{\infty}=D$. Moreover, it is clear that $D$, as an algebra, is generated by $K=D_0$ and $T=D_1/D_0$ with $D_1=K\oplus T$ if we identify an element $\xi={\xi}^id_i\in T$ with the vector field $\xi={\xi}^i(x){\partial}_i$ of differential geometry, but with ${\xi}^i\in K$ now. As a byproduct, 
the differential module $D^m$ is also filtred by the order and we obtain an {\it induced filtration} or {\it inductive limit} $0=M_{-1} \subseteq M_0 \subseteq M_1 \subseteq ... \subseteq M_q \subseteq ... \subseteq M_{\infty}=M$ with $d_iM_q\subseteq M_{q+1}$ provided by the prolongations. Now, {\it if we suppose that the system} $R_q=ker(\Phi)$ {\it is formally integrable}, then we have the {\it projective limit} $R=R_{\infty}\rightarrow ... \rightarrow R_q \rightarrow R_0\rightarrow 0$ obtained by successive jet projections. We have the crucial technical proposition ([28],[37],[43]):   \\

\noindent
{\bf PROPOSITION 2.A.2}: $R=hom_K(M,K)$ is a differential module for the Spencer operator and we have a bijective correspondence  $M_q \leftrightarrow  R_q=hom_K(M_q,K)$ over $K$ because $K$ is a field.  \\

\noindent
{\it Proof}: for any $f\in R$ and $m\in M$, we may set for any $f\in R, m\in M$:   \\
\[   (af)(m) =a(f(m))=f(am), \forall a\in K \hspace {1cm}  ,  \hspace{1cm} (\xi f)(m)=\xi (f(m)- f(\xi m), \forall \xi \in T   \]
and ckeck that we have successively with ${\xi}.a={\xi}^r{\partial}_ra$:   \\
\[   \begin{array}{rcl}
      ((\xi a)f)(m)  &  =   &  (\xi (af ))(m)   \\
                           & =  &   \xi ( af( m)) - af(\xi m)   \\
                           &  =  &  (\xi.a)f(m) + a(\xi.f(m)) -f(a\xi m)   \\                
                           &  =  &   (\xi.a)f(m)   + ((a\xi) f)(m)                                                
 \end{array}   \]                          
a result leading to $\xi a = a \xi +\xi.a$ in the operator sense. Setting finally $f(y^k_{\mu})=f^k_{\mu}$ with a slight abuse of notations when using the same notation $y^k_{\mu}$ for the residue instead of the standard ${\bar{y}}^k_{\mu}$. It follows that $R$ is a differential module for the law:  \\
\[  (d_if)(y^k_{\mu})= d_i(f(y^k_{\mu}) - f(d_iy^k_{\mu})={\partial}_if^k_{\mu}- f(y^k_{\mu+1_i})={\partial}_if^k_{\mu} - f^k_{\mu+1_i}   \]
and we have $d_id_jf=d_jd_if=d_{ij}f , \forall f\in R $ .  \\
\hspace*{12cm}  Q.E.D.   \\

   In this last section, we shall only deal wih linear or linearized differential operators. However, as explained with details in ([23],[27],[32],[34]), there is a nonlinear counterpart using the {\it nonlinear Janet sequence} coming from the {\it Vessiot structure equations} and a {\it nonlinear Spencer sequence}. However, the so-called {\it vertical machinery} involved, that is a systematic use of {\it fibered manifolds} and {\it vertical bundles}, is much more difficult though we have chosen the notations of this paper in such a way that the interested reader may easily adapt them. As for the Vessiot structure equations first found in 1903 ([46]), {\it they have been totally ignored during more than one century} for reasons that are not scientific at all (See [24] and the original letters presented in [25] for explanations).\\

Collecting all the results so far obtained, if a differential operator ${\cal{D}}$ is given in the framework of differential geometry, we may keep the same operator matrix in the framework of differential modules which are {\it left} modules over the ring $D$ of linear differential operators. We may also apply duality over $D$, that is apply $hom_D(\bullet,D)$, provided we deal now with {\it right} differential modules or use the operator matrix of $ad({\cal{D}})$ and deal again with {\it left} differential modules obtained through the $left \leftrightarrow right$ {\it conversion} procedure. In actual practice, it is essential to notice that {\it the new operator matrix may be quite different from the only transposed of the previous operator}, even if we are dealing with constant coefficients.\\

\noindent
{\bf DEFINITION 2.A.3}: If a differential operator $\xi \stackrel{\cal{D}}{\longrightarrow} \eta$ is given, a {\it direct problem} is to find (generating) {\it compatibility conditions} (CC) as an operator $\eta \stackrel{{\cal{D}}_1}{\longrightarrow} \zeta $ such that ${\cal{D}}\xi=\eta \Rightarrow {\cal{D}}_1\eta=0$. Conversely, given $\eta \stackrel{{\cal{D}}_1}{\longrightarrow} \zeta$, the {\it inverse problem} will be to look for $\xi \stackrel{\cal{D}}{\longrightarrow} \eta$ such that ${\cal{D}}_1$ generates the CC of ${\cal{D}}$ and we shall say that ${\cal{D}}_1$ {\it is parametrized by} ${\cal{D}}$ {\it if such an operator} ${\cal{D}}$ {\it is existing}.  \\

 \noindent
 {\bf REMARK 2.A.4}: Of course, solving the direct problem (Janet, Spencer) is {\it necessary} for solving the inverse problem. However, though the direct problem always has a solution, the inverse problem may not have a solution at all and the case of the Einstein operator is one of the best non-trivial PD counterexamples (Compare [27] to [47]). It is rather striking to discover that, in the case of OD operators, it took almost 50 years to understand that the possibility to solve the inverse problem was equivalent to the controllability of the corresponding control system (Compare [11] to [28]) and the situation is similar in GR as the above result has been first found in 1995 ([27],[34]).  \\
 
As $ad(ad(P))=P, \forall P \in D$, any operator is the adjoint of a certain operator and we recall that the {\it double duality test} needed in order to check whether $t(M)=0$ or not and to find out a parametrization if $t(M)=0$ when $M$ is defined by ${\cal{D}}_1$ has 5 steps which are drawn in the following diagram where $ad({\cal{D}})$ generates the CC of $ad({\cal{D}}_1)$ and ${\cal{D}}_1'$ generates the CC of 
${\cal{D}}=ad(ad({\cal{D}}))$:  \\
\[  \begin{array}{rcccccl}
 & & & & &  {\zeta}' &\hspace{15mm} 5  \\
 & & & & \stackrel{{\cal{D}}'_1}{\nearrow} &  &  \\
4 \hspace{15mm}& \xi  & \stackrel{{\cal{D}}}{\longrightarrow} &  \eta & \stackrel{{\cal{D}}_1}{\longrightarrow} & \zeta &\hspace{15mm}   1  \\
 &  &  &  &  &  &  \\
 &  &  &  &  &  &  \\
 3 \hspace{15mm}& \nu & \stackrel{ad({\cal{D}})}{\longleftarrow} & \mu & \stackrel{ad({\cal{D}}_1)}{\longleftarrow} & \lambda &\hspace{15mm} 2
  \end{array}  \]
\vspace*{3mm}

\noindent
{\bf THEOREM 2.A.5}: We have ${\cal{D}}_1$ parametrized by ${\cal{D}} \Leftrightarrow {\cal{D}}_1\simeq {\cal{D}}'_1 \Leftrightarrow t(M)=0 \Leftrightarrow ext^1(N)=0 $ in the differential module framework when $N$ is defined by $ad({\cal{D}}_1)$. These results do not depend on the finite free presentations of $M$ or $N$ (See [28] and [29] for more details). \\
 
\noindent
{\bf COROLLARY 2.A.6}: In the differential module framework, if $F_1 \stackrel{{\cal{D}}_1}{\longrightarrow} F_0 \stackrel{p}{\longrightarrow} M \rightarrow 0$ is a finite free presentation of $M=coker({\cal{D}}_1)$ and we already know that $t(M)=0$ by using the preceding test and Theorem, then we may obtain an exact sequence $F_1 \stackrel{{\cal{D}}_1}{\longrightarrow} F_0 \stackrel{{\cal{D}}}{\longrightarrow} E $ of free differential modules where ${\cal{D}}$ is the parametrizing operator, both with an inclusion $ M \subset E$ by chasing. However, there may exist other parametrizations $F_1 \stackrel{{\cal{D}}_1}{\longrightarrow} F_0 \stackrel{{\cal{D}}'}{\longrightarrow} E' $ called {\it minimal parametrizations} such that $coker({\cal{D}}')$ is a torsion module and we have thus $rk_D(M)=rk_D(E')$ (See [37],[41]).  \\

As shown by the next examples, the main difficulty met in OD or PD applications is that $ad({\cal{D}})$ may not be formally integrable at all, even if 
${\cal{D}}$ is involutive (See [31],[37-39] for other examples).  \\

\noindent
{\bf EXAMPLE 2.A.7}: Let us multiply on the left the second order trivially involutive single second order OD equation of the Introduction by a test function $\lambda$ and integrate by parts. The kernel of the adjoint operator is defined by:  \\
\[ y^2\,\,\, \rightarrow  \,\,\,d_{xx}\lambda - a(x)\lambda=0  \hspace{1cm}  , \hspace{1cm} y^3 \,\,\,\rightarrow \,\,\, - d_{xx}\lambda + d_x\lambda=0 \]
By addition we get $d_x\lambda  - a(x)\lambda=0$ and the adjoint OD system is not formally integrable. Differentiating once and substituting, we obtain the zero order OD equation $({\partial}_xa+a^2 - a)\lambda=0$ and the control system is controllable if and only if the ajoint operator is injective, that is if ${\partial}_xa +a^2-a\neq 0$.  \\

\noindent
{\bf EXAMPLE 2.A.8}:({\it Double pendulum}) If a rigid bar is able to move horizontally with reference position $x$ and we attach two pendula with respective length $l_1$ and $l_2$ making the (small) angles ${\theta}_1$ and ${\theta}_2$ with the vertical, the corresponding involutive control system is:  \\
\[d^2x +l_1d^2{\theta}^1 +g {\theta}^1=0, \hspace{1cm}  d^2x + l_2 d^2{\theta}^2 + g{\theta}^2=0   \]
where $g$ is the gravity. Multyplying these OD equations by two test functions ${\lambda}^1,{\lambda}^2$ and integrating by parts, we get the adjoint system:  \\
\[  x \rightarrow d^2{\lambda}^1 +d^2{\lambda}^2=0, \hspace{1cm}{\theta}^1 \rightarrow l_1d^2{\lambda}^1 + g{\lambda}^1=0, \hspace{1cm}  {\theta}^2 \rightarrow l_2d^2{\lambda}^2 + g{\lambda}^2=0\] 
Multiplying the secong equation by $l_2$, the third by $l_1$ while using the first, we obtain the zero order OD equation $l_2{\lambda}^1 + l_1{\lambda}^2=0$. Differentiating {\it twice} this time and substituting, we obtain the new zero order OD equation $(l_2/l_1){\lambda}^1 + (l_1/l_2) {\lambda}^2=0$. The determinant of the system of two zero order equations is then seen to be {\it exactly} $l_1 - l_2$. It follows that the system is controllable if and only if $l_1$ is different from $l_2$, a fact that the reader can check easily by himself when moving the bar conveniently. If one length depends on time, the corresponding controllability condition cannot be obtained without computer algebra, {\it even on such an elementary control system}. \\

\noindent
{\bf B) Variational calculus}  \\

\noindent
 {\bf EXAMPLE 2.B.1}: ({\it OD/PD Optimal Control Revisited}) Using the notations of the previous Formal Test, let us assume that the two differential sequences:  \\
\[  \begin{array}{rcccl}
 \xi & \stackrel{\cal{D}}{\longrightarrow} & \eta & \stackrel{{\cal{D}}_1}{\longrightarrow} & \zeta   \\
  \nu& \stackrel{ad(\cal{D})}{\longleftarrow} & \mu & \stackrel{ad({\cal{D}}_1)}{\longleftarrow} & \lambda 
  \end{array}  \] 
are {\it formally exact}, that is ${\cal{D}}_1$ generates the CC of ${\cal{D}}$ {\it and} $ad({\cal{D}})$ generates the CC of $ad({\cal{D}}_1)$, namely $\xi$ is a potential for ${\cal{D}}_1$ {\it and} $\lambda$ is a potential for $ad({\cal{D}})$. We may consider a variational problem for a cost function or lagrangian $\varphi (\eta)$ under the linear OD or PD constraint described by ${\cal{D}}_1\eta=0$. \\

\noindent
$\bullet$ Introducing convenient Lagrange multipliers $\lambda$ while setting $dx=dx^1\wedge ... \wedge dx^n$ for simplicity, we must vary the integral:  \\
 \[       \Phi=\int [\varphi(\eta) - \lambda {\cal{D}}_1\eta]dx  \Rightarrow \delta \Phi=\int [(\partial\varphi(\eta)/\partial\eta)\delta\eta - \lambda{\cal{D}}_1\delta\eta]dx \]
Integrating by parts, we obtain the EL equations: \\
\[  \partial\varphi(\eta)/\partial\eta = ad({\cal{D}}_1)\lambda  \]
to which we have to add the constraint ${\cal{D}}_1\eta=0$ obtained by varying $\lambda$ independently. If $ad({\cal{D}}_1)$ is an injective operator, in particular if ${\cal{D}}_1$ is formally surjective (no CC) while $n=1$ as inOD optimal control and $M$ is torsion-free, thus free ([13],[28]) or $n\geq 1$ and $M$ is projective, then one can obtain $\lambda$ {\it explicitly} and eliminate it by substitution. Otherwise, using the CC $ad({\cal{D}})$ of $ad({\cal{D}}_1)$, we have to study the formal integrability of the combined system: \\
\[     ad({\cal{D}})\partial\varphi(\eta)/\partial\eta=0, \hspace{4mm} {\cal{D}}_1\eta=0   \]
which may be a difficult task as we already saw through the examples of the Introduction.\\

\noindent
$\bullet$ However, {\it we may also transform the given variational problem with constraint to a variational problem without any constraint if and only if the differential constraint can be parametrized}. Using the parametrization of ${\cal{D}}_1$ by ${\cal{D}}$, we may vary the integral: \\
\[  \Phi=\int \varphi({\cal{D}}\xi)dx  \Rightarrow \delta \Phi = \int (\partial\varphi(\eta)/\partial\eta){\cal{D}}\delta\xi  dx                \]
whenever $\eta={\cal{D}}\xi$ and integrate by parts for arbitrary $\delta\xi$ in order to obtain the EL equations:  \\
\[   ad({\cal{D}})\partial\varphi(\eta)/\partial\eta=0, \hspace{4mm}  \eta={\cal{D}}\xi  \] 
 in a coherent way with the previous approach. \\
 
 As a byproduct, if the {\it field equations} ${\cal{D}}_1\eta=0$ can be parametrized by a {\it potential} $\xi$ through the formula ${\cal{D}}\xi=\eta$, then the {\it induction equations} $ad({\cal{D}})\mu=\nu$ can be obtained by duality in a coherent way with the {\it double duality test}, {\it on the condition to know what sequence must be used}. However, we have yet proved in ([31],[32],[35],[36],[39],[40]) that the {\it Cauchy stress equations} must be replaced by the {\it Cosserat couple-stress equations} and that the {\it Janet sequence} (only used in this paper) must be thus replaced by the {\it Spencer sequence}. Accordingly, it will become clear that the work of Lanczos ([14-17]) and followers ([1],[3-8],[18],[20-21],[42],[45]) using either ([2], exterior calculus), ([10], Janet and Gr\"{o}bner bases) or ([23], Pommaret bases) has been based on a {\it double confusion} between fields and inductions on one side, but also between the Janet sequence (only used in this paper) and the Spencer sequence. We have ([39], Proposition 5.1, p 2146): \\
 
\noindent
{\bf THEOREM  2.B.2}: The Janet and Spencer sequences for {\it any} Lie operator of finite type are formally exact by construction, {\it both with their corresponding adjoint sequences}. Lanczos has been trying to parametrize $ad({\cal{D}}_1)$ by $ad({\cal{D}}_2)$ when ${\cal{D}}_1$ is parametrizing ${\cal{D}}_2$. On the contrary, we have proved that one must parametrize $ad({\cal{D}})$ by $ad({\cal{D}}_1)$ when ${\cal{D}}$ is parametrizing ${\cal{D}}_1$ as in the famous {\it infinitesimal equivalence problem} ([23], p 332-336) or as in the above example. \\

\noindent
{\bf 3) APPLICATIONS}  \\

\noindent
{\bf EXAMPLE 3.1}: With $m=n=1, q=2$ we prove that the computation of the extension modules is also difficult for ordinary differential equations 
([34],[38] With $0\neq \alpha\in T^*$ and $\gamma$ transforming like the Christoffel symbols, we consider the second order geometric object $\omega=(\alpha,\gamma)$ and the second order system of infinitesimal Lie equations in Medolaghi form:  \\
\[  A\equiv \alpha {\partial}_x\xi + \xi {\partial}_x  \alpha=0, \,\,\, \Gamma\equiv {\partial}_{xx}\xi + \gamma {\partial}_x\xi + \xi {\partial}_x \gamma =0   \]
Multiplying $(A,\Gamma)$ by $({\mu}^1,{\mu}^2)$ and integrating by parts, we obtain $ad({\cal{D}})$ in the form:   \\
\[            {\partial}_{xx}{\mu}^2 -\alpha {\partial}_x{\mu}^1 - \gamma {\partial}_x{\mu}^2=\nu   \]
Now, we let the reader check by himself that $({\partial}_x\alpha- \alpha \gamma )/{\alpha}^2$ transforms like a scalar and the only Vessiot structure equation is ${\partial}_x\alpha  - \alpha\gamma  =\,\, c \,\, {\alpha}^2$ leading by linearization to the involutive system:   \\
\[                 {\partial}_x A - \gamma A - \alpha \Gamma - 2 \,\, c \,\, {\alpha} A=0   \]
and involutive operator ${\cal{D}}_1$. Multiplying on the left by the test function $\lambda$,  we obtain for $ad({\cal{D}}_1)$:   \\
\[  -{\partial}_x\lambda -(\gamma + 2 \, c \,\alpha) \lambda={\mu}^1, \,\,\, - \alpha \lambda ={\mu}^2  \]
Substituting $\lambda = - (1/\alpha){\mu}^2$, we discover that ${ext}^1(M)\neq 0$ is generated by the residue of the torsion element $(1/\alpha){\partial}_x{\mu}^2 +c{\mu}^2 - {\mu}^1={\nu}'$ satisfying ${\partial}_x{\nu}'=0$. Finally, as we have clearly $ker (ad({\cal{D}}_1))=0$, we obtain ${ext}^2(M)=0$.  \\

\noindent
{\bf EXAMPLE 3.2}: With $m=n=2, q=1, K=\mathbb{Q}<\omega>$ and $\omega=(\alpha,\beta)$ with $\alpha\in T^*,\beta\in {\wedge}^2T^*$, let us consider the Lie operator ${\cal{D}}:T \rightarrow  \Omega: \xi\rightarrow {\cal{L}} (\xi)\omega=(A={\cal{L}}(\xi)\alpha, B={\cal{L}}(\xi)\beta)$. The corresponding first order system:ÊÊ\\
\[ A_i \equiv {\alpha}_r{\partial}_i{\xi}^r+{\xi}^r{\partial}_r{\alpha}_i=0,  \,\,\, B \equiv \beta {\partial}_r{\xi}^r  +{\xi}^r{\partial}_r \beta =0\]       
is involutive whenever $\beta\neq 0$ and $d\alpha = c \beta$ where now $d$ is the standard exterior derivative and $c\in C$, exactly as in ([46], p 438-440) . We have the differential sequence:  \\
\[    0 \rightarrow \Theta \rightarrow T \stackrel{{\cal{D}}}{\longrightarrow}T^*{\times}_X{\wedge}^2T^* \stackrel{{\cal{D}}_1}{\longrightarrow} {\wedge}^2T^* \rightarrow 0  \]
or the resolution:  \\
\[    0 \rightarrow D \stackrel{{\cal{D}}_1}{\longrightarrow} D^3 \stackrel{{\cal{D}}}{\longrightarrow} D^3 \stackrel{p}{\longrightarrow} M  \rightarrow 0  \]
Multiplying $(A_1,A_2,B)$ respectively by $({\mu}^1,{\mu}^2,{\mu}^3)$, we obtain $ad({\cal{D}})$ in the form:  \\
\[ - {\alpha}_1({\partial}_1{\mu}^1+{\partial}_2{\mu}^2) - \beta ({\partial }_1{\mu}^3 - c {\mu}^2)={\nu}^1, \,\, 
   - {\alpha}_2({\partial}_1{\mu}^1+{\partial}_2{\mu}^2) - \beta ({\partial }_2{\mu}^3 + c {\mu}^1)={\nu}^2  \] 
Then, multiplying ${\partial}_1A_2-{\partial}_2A_1- cB$ by $\lambda$, we obtain $ad({\cal{D}}_1)$ as:  \\                       
\[   {\partial}_2\lambda={\mu}^1, \,\, - {\partial}_1 \lambda = {\mu}^2, \,\, - c \lambda={\mu}^3  \]
We have therefore to consider the two cases: \\
\noindent 
$\bullet$ $c=0$: We have the new CC ${\partial}_1{\mu}^1 + {\partial}_2{\mu}^2=0$ and ${\mu}^3=0$. It follows that the torsion module ${ext}^1(M \neq 0)$ is generated by the residue of ${\mu}^3={\nu}' $ because $\alpha \neq 0$ and we may thus suppose that ${\alpha}_1\neq 0$. As for ${ext}^2(M)$, this torsion module is just defined by the system ${\partial}_2\lambda=0, {\partial}_1\lambda=0$ for $\lambda$ and thus ${ext}^2(M)\neq 0$.  \\
\noindent
$\bullet$ $c\neq 0$: We must have the new CC:  \\
\[  {\partial}_1{\mu}^3 -  c{\mu}^2=0, {\partial}_2{\mu}^3 + c{\mu}^1=0\Rightarrow {\partial}_1{\mu}^1+{\partial}_2{\mu}^2=0  \]
It follows that ${ext}^1(M)$ is now generated by the residue of ${\partial}_1{\mu}^1 + {\partial}_2{\mu}^2={\nu}'$. Finally, $ker(ad({\cal{D}}_1))$ is defined by $\lambda=0$ and thus ${ext}^2(M)=0$. \\
Hence, both ${ext}^1(M) $ and ${ext}^2(M)$ highly depend on the Vessiot structure constant $c$.  \\

\noindent
{\bf EXAMPLE 3.3}: ({\it Contact transformations}) \\
 With $m=n=3, q=1, K=\mathbb{Q}(x^1,x^2,x^3)$  or simply $\mathbb{Q}(x)$, we may introduce the $1$-form $\alpha=dx^1-x^3dx^2 \in T^*$ and consider the system of finite Lie equations defined by $j_1(f)^{-1}(\alpha)=\rho(x) \alpha$. Eliminating the factor $\rho$ and linearizing at the $q$-jet of the idenity, we obtain the first order system of infinitesimal Lie equations:  \\
\[      {\Phi}^1\equiv {\partial}_2{\xi}^1 - x^3 {\partial}_2{\xi}^2 + x^3 {\partial}_1{\xi}^1 - (x^3)^2{\partial}_1{\xi}^2 - {\xi}^3=0, \,\,\, {\Phi}^2\equiv {\partial}_3{\xi}^1 - x^3 {\partial}_3{\xi}^2=0  \]
The first equation may provide ${\xi}^3$ whenever $({\xi}^1,{\xi}^2)$ are known from the second and there is therefore no CC for the correspoding Lie operator. However, this system is not even formally integrable because, differentiating the second equation with respect to $x^1$ and $x^2$ then substracting the  equation obtained by differentiating the first with respect to $x^3$, we obtain:   \\
\[      {\Phi}^3\equiv  x^3d_1{\Phi}^2 + d_2{\Phi}^2 - d_3{\Phi}^1\equiv  {\partial}_3{\xi}^3 + {\partial}_2{\xi}^2 + 2x^3{\partial}_1{\xi}^2=0  \]
which leads to an involutive system with two equations of class $3$, one equation of class $2$ and thus one CC of order $1$, namely $d_3{\Phi}^1-d_2{\Phi}^2 - x^3 d_1{\Phi}^2 + {\Phi}^3=0$. Now, it is well known that this {\it contact} operator ${\cal{D}}={\cal{D}}_0$ itself can be parametrized by an operator ${\cal{D}}_{-1}$  as follows:  \\
\[     - x^3 {\partial}_3\phi + \phi={\xi}^1, \,\,\,-{\partial}_3\phi={\xi}^2, \,\,\, {\partial}_2\phi+x^3{\partial}_1 \phi={\xi}^3 \,\,\, \Rightarrow \,\,\,  {\xi}^1 - x^3 {\xi}^2 = \phi\]
and thus $M\simeq D$. Accordingly, we have the two possible differential sequences:  \\
\[  \begin{array}{lcl}
0\rightarrow \Theta \rightarrow 3 \rightarrow 2  \rightarrow 0 & \Leftrightarrow & 0 \rightarrow D^2 \rightarrow D^3 \rightarrow M  \rightarrow 0   \\
 0 \rightarrow \Theta \rightarrow 3 \rightarrow 3 \rightarrow 1  \rightarrow 0  &  \Leftrightarrow  &  0 \rightarrow D \rightarrow D^3 \rightarrow D^3 \rightarrow M \rightarrow 0   
 \end{array}   \]
As $M$ is therefore free and thus projective, the two sequences of free modules split and, applying $hom_D(\bullet, D)$ we obtain split exact sequences of free right differential modules. Passing from right to left differential modules by the {\it side changing} procedure $N_D \rightarrow N={ }_DN= hom_K({\wedge}^nT^*,N_D)$. It follows that the adjoint sequence is exact too, though not strictly exact.  As such a result does not  depend on the differential sequence used, we prove it on the shortest sequence.  \\
For this, multiplying $({\Phi}^1,{\Phi}^2)$ by the test functions $({\mu}^1,{\mu}^2)$, contracting and integrating by parts, we obtain $ad({\cal{D}})$ in the form:  \\
\[  \left\{     \begin{array}{rrl}
{\xi}^1 \rightarrow & - {\partial}_2{\mu}^1 - x^3{\partial}_1{\mu}^1 - {\partial}_3{\mu}^2 &  = {\nu}^1  \\
{\xi}^2 \rightarrow &  x^3{\partial}_2{\mu}^1 +(x^3)^2{\partial}_1{\mu}^1 + x^3 {\partial}_3{\mu}^2 + 
{\mu}^2&  = {\nu}^2   \\
 {\xi}^3 \rightarrow  &  -  {\mu}^1  & ={\nu}^3
 \end{array}   \right.  \]
with ${\nu}^2+x^3{\nu}^1={\mu}^2$. We obtain therefore $ker(ad({\cal{D}}))=0 \Rightarrow {ext}^1(M)=0$ and we have trivially ${ext}^2(M)={ext}^3(M)=0$.   \\
Coming back to the Vessiot structure equations, we notice that $\alpha$ is not invariant by the contact Lie pseudogroup and cannot be considered as an associated geometric object. We have shown in ([24], p 684-691) that the corresponding geometric object is a $1$-form density $\omega$ leading to the system of infinitesimal Lie equations in Medolaghi form:  \\
\[            {\Omega}_i\equiv ( {\cal{L}}(\xi)\omega)_i \equiv {\omega}_r{\partial}_i{\xi}^r - \frac{1}{2}{\omega}_i{\partial}_r{\xi}^r + {\xi}^r{\partial}_r{\omega}_i=0  \]
and to the {\it only} Vessiot structure equations:  \\
\[        {\omega}_1({\partial}_2{\omega}_3- {\partial}_3{\omega}_2) +   {\omega}_2({\partial}_3{\omega}_1- {\partial}_1{\omega}_3) +  {\omega}_3({\partial}_1{\omega}_2- {\partial}_2{\omega}_3)=c  \]
with the {\it only} structure constant $c$. In the present contact situation, we may choose 
$\omega=(1, -x^3,0)$ and get $c=1$ but we may also choose $\omega=(1,0,0)$ and get $c=0$, these two choices both bringing an involutive system. Let us prove that the situation becomes completely different with the new system:  \\
\[  -2{\Omega}_1\equiv  {\partial}_3{\xi}^3 + {\partial}_2{\xi}^2 - {\partial}_1{\xi}^1=0,\,\, {\Omega}_2\equiv {\partial}_2{\xi}^1=0, \,\, {\Omega}_3\equiv {\partial}_3{\xi}^1=0   \]
having the only CC $ d_2{\Omega}_3 - d_3 {\Omega}_2=0 $. \\
Multilying the three previous equations by the three test functions $\mu$, the only CC by the test function $\lambda$ and integrating by parts, we get the adjoint operators:    \\
\[    0= {\mu}^1, \,\,  {\partial}_3\lambda= {\mu}^2, \,\, - {\partial}_2\lambda={\mu}^3   \]
\[ {\partial}_1{\mu}^1 - {\partial}_2{\mu}^2 - {\partial}_3{\mu}^3= {\nu}^1, \,\, - {\partial}_2{\mu}^1={\nu}^2, \,\, - {\partial}_3{\mu}^1={\nu}^3 \]
It follows that $0\neq D{\xi}^1 =t(M) \subset M$ with a strict inclusion and ${ext}^1(M)\neq 0$. Similarly, $ker(ad({\cal{D}}_1))$ is defined by ${\partial}_2\lambda=0, {\partial}_3\lambda=0$ and thus ${ext}^2(M)\neq 0$.  \\
Our problem will be now to construct and compare the differential sequences:   \\
\[   \phi \stackrel{{\cal{D}}_{-1}}{\longrightarrow} \xi \stackrel{{\cal{D}}}{\longrightarrow} \Omega \stackrel{{\cal{D}}_1}{\longrightarrow}C \]  
\[   \theta \stackrel{ad({\cal{D}}_{-1})}{\longleftarrow} \nu \stackrel{ad({\cal{D}})}{\longleftarrow} \mu \stackrel{ad({\cal{D}}_1)}{\longleftarrow} \lambda  \]                                                                                                                                                                                                                                                      
For this, linearizing the only Vessiot structure equation, we get the CC operator ${\cal{D}}_1$ and the corresponding system ${\cal{D}}_1 \Omega=0$ in the form:  \\
\[   {\omega}_1({\partial}_2{\Omega}_3- {\partial}_3{\Omega}_2) +   {\omega}_2({\partial}_3{\Omega}_1- {\partial}_1{\Omega}_3) +  {\omega}_3({\partial}_1{\Omega}_2  - {\partial}_2{\Omega}_1) \]
\[+ ({\partial}_2{\omega}_3- {\partial}_3{\omega}_2) {\Omega}_1
  + ({\partial}_3{\omega}_1- {\partial}_1{\omega}_3) {\Omega}_2+  ({\partial}_1{\omega}_2- {\partial}_2{\omega}_3) {\Omega}_3=0   \]
Multiplying on the left by the test function $\lambda$ and integrating by parts, we get the operator $ad({\cal{D}}_1)$ in the form:  \\
\[  \left\{   \begin{array}{lccl}
{\Omega}_1 \rightarrow  & {\omega}_3{\partial}_2\lambda - {\omega}_2 {\partial}_3\lambda + 2({\partial}_2{\omega}_3 - {\partial}_3 {\omega}_2)\lambda & = {\mu}^1  \\
{\Omega}_2 \rightarrow  & {\omega}_1{\partial}_3\lambda - {\omega}_3 {\partial}_1\lambda + 2({\partial}_3{\omega}_1 - {\partial}_1 {\omega}_3)\lambda & = {\mu}^2  \\
{\Omega}_3 \rightarrow  & {\omega}_2{\partial}_1\lambda - {\omega}_1 {\partial}_2\lambda + 2({\partial}_1{\omega}_2 - {\partial}_2 {\omega}_1)\lambda & = {\mu}^3  \\
\end{array} \right.  \]
We obtain therefore the crucial formula $2c \,\lambda = {\omega}_i{\mu}^i$ showing how the previous sequences are {\it essentially} depending on the Vessiot structure constant $c$. Indeed, if $c\neq 0$, then $\mu=0 \Rightarrow \lambda =0$ and the operator $ad({\cal{D}}_1)$ is injective. This is the case when $\omega =(1,-x^3,0)\Rightarrow c=1  \Rightarrow \lambda =0$. On the contrary, if $c=0$, then the operator $ad({\cal{D}}_1)$ may not be injective as can be seen by choosing $\omega=(1,0,0)$. Indeed, in this case we get a kernel defined by ${\partial}_3\lambda=0, {\partial}_2\lambda=0$.  \\
Finally, in order to exhibit the generating CC of $ad({\cal{D}}_1)$ when $c\neq 0$, we just need to substitute $\lambda=(1/2c){\omega}_i{\mu}^i$ in the previous equations $ad({\cal{D}}_1)\lambda=\mu$. On the other side, multiplying the equations ${\cal{D}}\xi=\Omega$ by test functions ${\mu}^i$ and integrating by parts, we get $ad({\cal{D}})\mu=\nu$ in the form:  \\
\[ \left\{  \begin{array}{lccl}
{\xi}^1 \rightarrow & \,\,  & - {\partial}_i({\omega}_1{\mu}^i) + \frac{1}{2}{\partial}_1({\omega}_i{\mu}^i)+
({\partial}_1{\omega}_i){\mu}^i&={\nu}^1    \\
{\xi}^1 \rightarrow & \,\,  & - {\partial}_i({\omega}_2{\mu}^i) + \frac{1}{2}{\partial}_2({\omega}_i{\mu}^i)+
({\partial}_2{\omega}_i){\mu}^i & = {\nu}^2   \\
{\xi}^1 \rightarrow & \,\,  & - {\partial}_i({\omega}_3{\mu}^i) + \frac{1}{2}{\partial}_3({\omega}_i{\mu}^i)+
({\partial}_3{\omega}_i){\mu}^i & = {\nu}^3
\end{array}  \right. \]
We let the reader check directly, as a delicate exercise by hand or with computer algebra, that we have indeed $ad({\cal{D}})\circ ad({\cal{D}}_1)\equiv 0$ and it remains to prove that $ad({\cal{D}})$ generates the CC of $ad({\cal{D}}_1)$. It is in such a situation that we can measure the usefulness of homological algebra and we only prove this result directly when $\omega=(1,-x^3,0) \Rightarrow c=1$. In this case, the kernel of $ad({\cal{D}})$ is easily seen to be defined by:  \\
\[  \begin{array}{ll}
  \frac{1}{2}{\partial}_1{\mu}^1 +{\partial}_2{\mu}^2 + {\partial}_3{\mu}^3+\frac{1}{2} x^3{\partial}_1{\mu}^2 & =0  \\
  - (x^3)^2{\partial}_1{\mu}^2 + x^3{\partial}_1{\mu}^1+{\partial}_2{\mu}^1-x^3{\partial}_2{\mu}^2+2{\mu}^3 & =0   \\
  {\partial}_3{\mu}^1 - x^3 {\partial}_3{\mu}^2 - 3 {\mu}^2  &  =0
  \end{array}  \]
while, setting now $2\lambda={\mu}^1-x^3{\mu}^2$ and substituting, the CC of $ad({\cal{D}}_1)$ seem to be only defined by the two PD equations:  \\
\[ \begin{array}{ll}
{\partial}_3{\mu}^1 - x^3 {\partial}_3{\mu}^2 - 3 {\mu}^2   &  =0  \\
{\partial}_2{\mu}^1 - x^3 {\partial}_2{\mu}^2 - (x^3)^2 {\partial}_1{\mu}^2 +x^3 {\partial}_1{\mu}^1+ 2{\mu}^3  &  =0
\end{array}  \]
The strange fact is that {\it such a system is not formally integrable} and one has to differentiate the second PD equation with respect to $x^3$ and substract the first PD equation differentiated with respect to $x^2$ in order to get the additional PD equation:   \\
\[       {\partial}_3{\mu}^3 +{\partial}_2{\mu}^2 + \frac{1}{2}x^3 {\partial}_1{\mu}^2+\frac{1}{2}{\partial}_1{\mu}^1=0     \]
in order to find an isomorphic involutive system, a result showing that the differential sequence and its formal adjoint are both formally exact though {\it not} strictly exact. We conclude this example with the following striking result:  \\

\noindent
{\bf THEOREM 3.4}: The {\it contact differential sequence} and its formal adjoint are both split long exact sequences of free and thus projective modules, if and only if $c\neq 0$. Moreover, the central operator ${\cal{D}}$ is formally adjoint with a slight abuse of language in the sense that a linear change of bases must be done.  \\  

\noindent
{\it Proof}: As we have just proved that $ad({\cal{D}})$ was generating the CC of $ad({\cal{D}}_1)$, we may look for the CC of $ad({\cal{D}})$ in order to recover the parametrization of ${\cal{D}}$ given at the beginning of this example.

We may write the operator ${\cal{D}} \xi ={\cal{L}}(\xi)\omega=\Omega$ in the form:  \\

\[  {\omega}_r{\partial}_i{\xi}^r - \frac{1}{2}{\omega}_i {\partial}_r{\xi}^r +{\xi}^r{\partial}_r{\omega}_i= {\Omega}_i  \]
that may be also written as:  \\
\[{\partial}_i({\omega}_r{\xi}^r) - \frac{1}{2} {\omega}_i{\partial}_r{\xi}^r + {\xi}^r({\partial}_r{\omega}_i - {\partial}_i{\omega}_r)={\Omega}_i\]

Linearising the Vessiot structure equation over $\omega$, we may write the CC operator ${\cal{D}}_1{\Omega}= C$ as:   \\
\[ \begin{array}{lcr}
 {\omega}_1({\partial}_2{\Omega}_3 - {\partial}_3{\Omega}_2) + {\omega}_2({\partial}_3{\Omega}_1 - {\partial}_1{\Omega}_3) + {\omega}_3({\partial}_1{\Omega}_2 - {\partial}_2{\Omega}_1) &  &   \\
+ ({\partial}_2{\omega}_3 - {\partial}_3{\omega}_2){\Omega}_1 + ({\partial}_3{\omega}_1 - {\partial}_1{\omega}_3){\Omega}_2 + ({\partial}_1{\omega}_2 - {\partial}_2{\omega}_1){\Omega}_3 &= &C 
 \end{array} \]

Our aim is now to study the differential sequence where we could set ${\cal{D}}={\cal{D}}_0$  in the sequences:  \\
\[   0 \rightarrow \phi \stackrel{{\cal{D}}_{-1}}{\longrightarrow}  \xi  \stackrel{{\cal{D}}}{\longrightarrow} \Omega \stackrel{{\cal{D}}_1}{\longrightarrow} C  \rightarrow 0   \]
and its formal dual:  \\
\[  0 \leftarrow \theta \stackrel{{ad(\cal{D}}_{-1})}{\longleftarrow} \nu \stackrel{ad(\cal{D})}{\longleftarrow} \mu \stackrel{ad({\cal{D}}_1)}{\longleftarrow} \lambda  \leftarrow 0  \]
First of all, as already noticed, we have ${\omega}_i{\mu}^i=2c\lambda$ and $ad({\cal{D}}_1)$ is injective if and only if $c\neq 0$. It follows that the differential module defined by ${\cal{D}}_1$ is projective and the sequences split if we are able to construct ${\cal{D}}_{-1}$ and to prove that it is an injective operator. For this, we notice that the symbol map of $ad({\cal{D}})$ is:  \\
\[    - {\omega}_i{\mu}^r_r + \frac{1}{2}{\omega}_r{\mu}^r_i={\nu}^i    \]
Hence, after one prolongation on the symbol level, we get the only CC:  \\
\[  {\omega}_1 ({\nu}^2_3 - {\nu}^3_2) + {\omega}_2 ({\nu}^3_1 - {\nu}^1_3) + {\omega}_3 ( {\nu}^1_2 - {\nu}^2_1)=0  \]
After substitution and a painfull computation or the use of computer algebra, one finally obtains:  \\
\[   \begin{array}{rcl}
{\omega}_1 ({\nu}^2_3 - {\nu}^3_2) + {\omega}_2 ({\nu}^3_1 - {\nu}^1_3) + {\omega}_3 ( {\nu}^1_2 - {\nu}^2_1)&  & \\
+2({\partial}_2{\omega}_3 - {\partial}_3 {\omega}_2) {\nu}^1 + 2({\partial}_3{\omega}_1 - {\partial}_1 {\omega}_3) {\nu}^2 + 2({\partial}_1{\omega}_2 - {\partial}_2{\omega}_1) {\nu}^3 & = & \theta    
\end{array}    \]
Wefinally obtain ${\cal{D}}_{-1}$ in the form:  \\
\[ \left\{    \begin{array}{lccl}
 {\omega}_2{\partial}_3\phi- {\omega}_3 {\partial}_2\phi + ({\partial}_2{\omega}_3 - {\partial}_3 {\omega}_2)\phi & = {\xi}^1  \\
 {\omega}_3{\partial}_1\phi - {\omega}_1 {\partial}_3\phi + ({\partial}_3{\omega}_1 - {\partial}_1 {\omega}_3)\phi & = {\xi}^2  \\
 {\omega}_1{\partial}_2\phi - {\omega}_2 {\partial}_1\phi + ({\partial}_1{\omega}_2 - {\partial}_2 {\omega}_1)\phi & = {\xi}^3  \\
\end{array}  \right. \]
and this operator is injective whenever $c\neq 0$ because ${\omega}_r{\xi}^r=c \phi$.  \\
The situation is similar in arbitrary dimension $n=2p+1$ with $1$-form $\alpha= dx^n -{\sum}^p_{\alpha=1}x^{\bar{\alpha}}dx^{\alpha}$ as we 
have again one Vessiot structure constant $c$ and the injective parametrization :  \\
\[  \phi - x^{\bar{\beta}}\frac{\partial \phi}{\partial x^{\beta}}= {\xi}^n, \hspace{6mm}  - \frac{\partial \phi}{\partial x^{\bar{\alpha}}}={\xi}^{\alpha}, \hspace{6mm}  \frac{\partial \phi}{\partial x^{\alpha}}+x^{\bar{\alpha}} \frac{\partial \phi}{\partial x^n}={\xi}^{\bar{\alpha}} \hspace{4mm}\Rightarrow \hspace{4mm}  \phi=\alpha(\xi)  \]
We have the locally exact Janet type split sequence:  \\
\[  0  \rightarrow {\wedge}^0T^* \stackrel{{\cal{C}}}{\longrightarrow} T \stackrel{\cal{D}}{\longrightarrow} F_0 \stackrel{{\cal{D}}_1}{\longrightarrow} F_1 \stackrel{{\cal{D}}_2}{\longrightarrow} ..... \stackrel{{\cal{D}}_{n-2}}{\longrightarrow} F_{n-2} \rightarrow 0  \]
with $dim(F_r)=n!/(r+2)!(n-r-2)!$ and we refer again the reader to ([24]) for more details because, when $n\geq  5$, we have to use a $2$-contravariant skewsymmetric tensor density. \\

\noindent
{\bf EXAMPLE 3.5}: ({\it Unimodular contact transformations })      \\
With $m=n=3, q=1$ and $K=\mathbb{Q}(x^1,x^2,x^3)$, let us consider the Lie pseudogroup of transformations $y=f(x)$ preserving the $1$-form $\alpha = dx^1 - x^3 dx^2$. It is defined by the Paffian system $dy^1 - y^3dy^2 = dx^1 - x^3d x^2$. After linearization over the identity transformation $y=x$, the corresponding Lie operator is ${\cal{D}}:\xi \rightarrow {\cal{L}}(\xi)\alpha =A$ by introducing the standard Lie derivatives on exterior forms. With $A=(u,v,w)$, the corresponding linear inhomogeneous system can be written under the form:\\

\[  \left\{  \begin{array}{lcl}
{\xi}^1_3 - x^3 {\xi}^2_3 & = & u \\
{\xi}^1_2 - x^3 {\xi}^2_2 -{\xi}^3 & = & v \\
{\xi}^1_1 - x^3 {\xi}^2_1 & = &  w 
\end{array}
\right. \fbox{$\begin{array}{lll}
1 & 2 & 3 \\
1 & 2 & \bullet \\
1 & \bullet & \bullet
\end{array}$}  \]
with one equation of class $3$, one equation of class $2$ and one equation of class $1$. Using the {\it Janet board} of {\it multiplicative /nonmultiplicative} variables and the corresponding so-called Pommaret basis ([30]), we discover at once that such a system is not involutive and thus not formally integrable as we may add $3$ new first order equations. Accordingly, ${\xi}^2$ and ${\xi}^3$ cannot be given arbitrarily when the right members vanish and {\it it is not evident at all} to discover that these three PD equations are differentially independent in the sense that there is no CC. Permuting the order of the independent variables so that $x^2 < x^3<x^1$ or permuting the independent coordinates with $(1,2,3)\Rightarrow (3,1,2)$, we obtain the following involutive system in solved form:  \\

\[  \left\{  \begin{array}{lcl}
{\xi}^2_1 & = & u_1- w_3  \\
{\xi}^3_1  & = &  w_2-v_1  \\
{\xi}^1_1 & = &  w - x^3 w_3 + x^3 u_1  \\
{\xi}^3_3 + {\xi}^2_2  & = & u_2 - v_3  \\
{\xi}^1_3 - x^3 {\xi}^2_3 & = & u \\
{\xi}^1_2 - x^3 {\xi}^2_2 -{\xi}^3 & = & v \\
\end{array}
\right. \fbox{$\begin{array}{lll}
2 & 3 & 1 \\
2 & 3 & 1  \\
2 & 3 & 1 \\
2 & 3 & \bullet  \\
2 & 3  & \bullet   \\
2 & \bullet & \bullet
\end{array}$}  \]

From the general formal theory of systems of PD equations, we know that the CC for $(u,v,w)$ can be deduced at once from the four first order CC obtained from the Janet board. Surprisingly, they do not produce any CC between $(u,v,w)$ as we can check:   \\

\[  d_1 (u) - d_3 (w - x^3 w_3 + x^3 u_1) + x^3 d_3 (u_1 - w_3) \equiv 0  \]
\[  .....................................................................\]
\[ d_3 (v)  + d_2 (u)  + (u_2 - v_3)  \equiv  0  \]

Finally, we notice that both $\{{\xi}^1,{\xi}^2, {\xi}^3\}$ are torsion elements and the differential module $M$ defined by ${\cal{D}}$ is a torsion module. \\
We now turn to the study of the Vessiot structure equations with two structure constants ([23-26],[38]). For this, if $0\neq \alpha \in T^*$ and $\beta \in {\wedge}^2T^*$ are such that $0\neq  \alpha \wedge \beta \in {\wedge}^3T^*$, we may introduce the geometric object $\omega=(\alpha, \beta)$ and the Vessiot structure equations become:  \\
\[      d\alpha = c' \beta, \hspace{2cm}  d\beta = c" \alpha \wedge \beta \hspace{1cm}  \Rightarrow \hspace{1cm}   c'c"=0  \]
by closing the exterior system. Accordingly, we must distinguish two cases:   \\

\noindent
$\bullet \hspace{3mm}c'=0$: Introducing the interior multiplication $i( \hspace{2mm}  )$ on forms, it is well known that:  \\
\[  d \alpha = 0 \hspace{4mm}\Rightarrow \hspace{3mm}  {\cal{L}}(\xi)\alpha=i(\xi)d\alpha + d i(\xi) \alpha=d(\alpha(\xi))=0  \]
and it follows that $\alpha(\xi)$ is a torsion element, namely $(1/x^1){\xi}^1$ if we set $\omega =(\alpha=(1/x^1)dx^1, \beta=x^1 dx^2\wedge dx^3) \Rightarrow c'=0, c"=1$ or simply ${\xi}^1$ if we set $\omega=(\alpha=dx^1, \beta=dx^2\wedge dx^3) \Rightarrow c'=0, c"=0$. It is important to notice that, in both cases, the differential module $M$ is defined not only by ${\cal{L}}(\xi)\alpha=0$, but also by ${\cal{L}}(\xi)\beta=0$ in such a way to have the finite free resolution:  \\
\[    0 \rightarrow D \stackrel{{\cal{D}}_2}{ \longrightarrow} D^4 \stackrel{{\cal{D}}_1}{\longrightarrow} D^6 \stackrel{{\cal{D}}}{\longrightarrow} D^3 \rightarrow M  \rightarrow 0   \]
in such a way that $rk_D(M)=3 - 6 - 4 - 1=0 $ and $M$ is thus a torsion module. Considering in particular the last case, we may easily define ${\cal{D}}$ by:  \\
\[ d_1{\xi}^3={\eta}^1, d_1{\xi}^2={\eta}^2, d_1{\xi}^1={\eta}^3, d_2{\xi}^2+d_3{\xi}^3={\eta}^4, d_2{\xi}^1={\eta}^5, d_3{\xi}^1={\eta}^6  \]
and ${\cal{D}}_1$ by:  \\
\[  d_1{\eta}^4-d_2{\eta}^2- d_3{\eta}^1={\zeta}^1, d_1{\eta}^5-d_2{\eta}^3={\zeta}^2, d_1{\eta}^6 - d_3{\eta}^3={\zeta}^3, d_2{\eta}^6 - d_3{\eta}^5={\zeta}^4  \]
in such a way that ${\cal{D}}_2$ is defined by $d_1{\zeta}^4-d_2{\zeta}^3+d_3{\zeta}^2=0$. We may then define $ad({\cal{D}})$ (up to sign) by: \\
\[  d_1{\mu}^3+d_2{\mu}^5+d_3{\mu}^6={\nu}^1, d_1{\mu}^2+d_2{\mu}^4={\nu}^2, d_1{\mu}^1+d_3{\mu}^4={\nu}^3  \] 
then $ad({\cal{D}}_1)$ by the Pommaret basis:  \\
\[   \left\{    \begin{array}{lcl}
 d_3{\lambda}^1& = &{\mu}^1 \\
 d_3{\lambda}^3+d_2{\lambda}^2 & = &  {\mu}^3 \\
 d_3{\lambda}^4 - d_1 {\lambda}^2& = &  {\mu}^5 \\
 d_2{\lambda}^1& = &{\mu}^2 \\
 d_2{\lambda}^4 + d_1{\lambda}^3& = &  - {\mu}^6 \\
d_1{\lambda}^1 & = & -{\mu}^4 
\end{array}
\right. \fbox{$\begin{array}{lll}
1& 2 & 3 \\
1 & 2 & 3  \\
1 & 2 & 3 \\
1 & 2 & \bullet  \\
1 & 2  & \bullet   \\
1 & \bullet & \bullet
\end{array}$}  \]
This involutive system has $4$ CC, namely the $3$ CC described by $ad({\cal{D}})$ of course, plus the only additional CC $\tau\equiv d_3{\mu}^2 - d_2{\mu}^1=0$ which brings the {\it only} torsion element $\tau$ because we have $d_1\tau= d_{13}{\mu}^2 - d_{12}{\mu}^1= d_3(d_1{\mu}^2+d_2{\mu}^4) - d_2(d_1{\mu}^1+d_3{\mu}^4)=0$ by taking the residue. It follows that he torsion module $ext^1(M)$ has the {\it only} generator $\tau$. \\
Finally, $ad({\cal{D}}_2)$ is defined by:  \\
\[     0={\lambda}^1, - d_3 \theta={\lambda}^2, d_2 \theta= {\lambda}^3, - d_1 {\theta}= {\lambda}^4  \]
and the torsion module ${ext}^2(M)\neq 0$ is generated by the residue of ${\lambda}^1$.\\

\noindent
$\bullet \hspace{3mm}c'\neq 0$: In this case, {\it we need to have} $c"=0$ and may suppose that $c'=1 \Rightarrow d\alpha=\beta \Rightarrow d\beta=0$ as before with $\omega=(\alpha=dx^1-x^3dx^2, \beta=dx^2\wedge dx^3)$. Also, as  $\alpha \neq 0$, we may also suppose that ${\alpha}_1\neq 0$. We may therefore only consider the generating equations ${\alpha}_r{\partial}_i {\xi}^r+{\xi}^r {\partial}_r{\alpha}_i=0, \forall i=1,2,3$. Multiplying by test functions $({\lambda}^1,{\lambda}^2, { \lambda}^3)$ and integrating by parts as usual, we get the adjoint operator ${\xi}^r \rightarrow  - {\alpha}_r({\partial}_i{\lambda}^i )+{\beta}_{ri}{\lambda}^i  =  {\mu}^r $. Looking at its kernel, we have ${\partial}_i{\lambda}^i={\beta}_{1i}{\lambda}^i$ that we may substitute in the two other PD equations, obtaining therefore $2$ linear equations that could be used in order to express $({\lambda}^2,{\lambda}^3) $in terms of ${\lambda}^1$ only, that we can substitute again in the first PD equations. This result proves that the differential module $N$ defined by $ad({\cal{D}})$ cannot vanish and is a torsion module with at least one generator, namely the residue of ${\lambda}^1$. It follows that $M$ cannot be free or even projective ([13], p 212). Nevertheless, we have $hom_D(M,D)=0$ and we have the dual short exact sequences:  \\
\[    0  \rightarrow D^3 \stackrel{{\cal{D}}}{\longrightarrow} D^3 \rightarrow M \rightarrow 0 , \hspace{10mm}
0 \leftarrow N \leftarrow D^3 \stackrel{ad({\cal{D}})}{\longleftarrow} D^3 \leftarrow 0    \]
even if we already know that ${\cal{D}}$ defined only by ${\cal{L}}(\xi)\alpha=0$ is far from being involutive as only the system defined by ${\cal{L}}(\xi)\alpha=0, {\cal{L}}(\xi)\beta=0$ is indeed involutive. Hence, both $M$ and $N$ are torsion modules with $3$ generators in such a way that ${ext}^1(M)=N\neq 0$ and ${ext}^1(N)=M \neq 0$. Finally, as $M$ admits a resolution with only one operator, then ${ext}^2(M)=0$ and it follows from this example that the structure of the ${ext}^i(M)$ highly depends on the Vessiot structure constants.  \\

\noindent
{\bf EXAMPLE 3.6}: ({\it Riemann/Lanczos problem})  \\
We now consider with details the Riemann/Lanczos problem which is at the same time the simplest of the two Lanczos problems as it can be solved in arbitrary dimension $n\geq 2$ but is also an example of seven successive confusions that have been done during the last fifty years. According to the last Section, the starting motivating point seems absolutely natural at first. Indeed, considering the {\it Killing} operator ${\cal{D}}:\xi \rightarrow {\cal{L}}(\xi)\omega=\Omega \in S_2T^*=F_0$ where ${\cal{L}}(\xi)$ is the Lie derivative with respect to $\xi$ and $\omega \in S_2T^*$ is a nondegenerate metric with $det(\omega)\neq 0$. Accordingly, it is a lie operator with ${\cal{D}}\xi=0, {\cal{D}}\eta=0 \Rightarrow {\cal{D}}[\xi,\eta]=0$ and we denote simply by $\Theta \subset T$ the set of solutions with $[\Theta,\Theta ]\subset \Theta$. Now, as we have explained many times, the main problem is to describe the CC of ${\cal{D}}\xi=\Omega \in F_0$ in the form ${\cal{D}}_1\Omega=0$ by introducing the so-called {\it Riemann} operator ${\cal{D}}_1:F_0 \rightarrow F_1$, using the standard notations that can be found at length in our many books ([23-26],[38]) or papers ([34],[37],[39],[40]). We advise the reader to follow closely the next lines and to imagine why it will not be possible to repeat them for studying the Weyl/Lanczos problem. Introducing the well known Levi-Civita isomorphism $j_1(\omega)=(\omega, {\partial}_x \omega) \simeq (\omega , \gamma)$ by 
introducing the Christoffel symbols ${\gamma}^k_{ij}=\frac{1}{2}{\omega}^{kr}({\partial}_i{\omega}_{rj} + {\partial}_j{\omega}_{ir}-
{\partial}_r{\omega}_{ij})$ where $({\omega}^{rs})$ is the inverse matrix of $({\omega}_{ij})$, we get $R_2 \subset J_2(T)$:  \\
\[  \left\{  \begin{array}{lcl}
{\Omega}_{ij}  & \equiv & {\omega}_{rj}(x){\xi}^r_i + {\omega}_{ir}(x){\xi}^r_j + {\xi}^r {\partial}_r {\omega}_{ij}(x) =0  \\
{\Gamma}^k_{ij} & \equiv  & {\xi}^k_{ij} + {\gamma}^k_{rj}(x) {\xi}^r_i+{\gamma}^k_{ir}(x) {\xi}^r_j +{\gamma}^k_{ir}(x) {\xi}^r_j - {\gamma}^r_{ij}(x) {\xi}^k_r + {\xi}^r{\partial}_r {\gamma }^k_{ij}(x)=0
\end{array} \right.  \]
if we use jet coordinteswith sections ${\xi}_q:x \rightarrow ({\xi}^k(x), {\xi}^k_i(x), {\xi}^k_{ij}(x), ... )$ transforming like $j_q(\xi): x \rightarrow ({\xi}^k(x), {\partial}_i{\xi}^k(x), {\partial}_{ij}{\xi}^k(x), ...)$. The system $R_1 \subset J_1(T)$ has a symbol $g_1 \simeq {\wedge}^2T^*\subset T^* \otimes T$ depending only on $\omega$ with $dim(g_1)=n(n-1)/2$ and is finite type because its first prolongation is $g_2=0$. It cannot be thus involutive as can be seen directly on the following Janet board for finding a Pommaret basis when $n=2$ and $\omega$ is the euclidean metric:  \\
\[  \left\{  \begin{array}{l}
{\xi}^2_2=0  \\
{\xi}^1_2 + {\xi}^2_1=0  \\
{\xi}^1_1=0
\end{array} \right. \fbox{$\begin{array}{ll}
1 & 2  \\
1 & 2   \\
1 & \bullet 
\end{array}$}  \]
Indeed, the only dot appearing in the board cannot provide any CC for the symbol $g_1$ and we have therefore the short exact sequence:  \\
\[   0 \rightarrow g_2 \rightarrow S_2T^*\otimes T \rightarrow T^* \otimes F_0 \rightarrow 0  \]
by using the fact that $g_2=0$ and counting the common dimension $n(n+1)/2$, because an epimorphism between two spaces of the same dimensionis also a monomorphism and thus an isomorphism. Accordingly, we need to use one additional prolongation and arrive to the:  \\

\noindent
$\bullet$ {\it First confusion}: Using one of the main results to be found in ([37,[38],....), we know that, when $R_1$ is formally integrable, then the CC of ${\cal{D}}$ are of order $s+1$ where $s$ is the number of prolongations needed in order to get an involutive symbol, that is $s=1$ in the present situation, a result that should lead to CC of order $2$ if $R_1$ were formally integrable.

\[  \begin{array}{rcccccccl}
& & & 0 && 0 & & &  \\
& & & \downarrow & & \downarrow & & &  \\
 & 0 & \rightarrow & S_3T^* \otimes T & \rightarrow & S_2T ^*\otimes F_0  & \rightarrow & h_2  &\rightarrow 0  \\
 &  \downarrow &  &  \downarrow & & \downarrow & & \downarrow &  \\
0 \rightarrow  & R_3 & \rightarrow  & J_3(T) & \rightarrow & J_2(F_0)&  \rightarrow & Q_2 & \rightarrow 0  \\
  &  \downarrow &  &  \downarrow & & \downarrow & & \downarrow &  \\
0 \rightarrow  & R_2 & \rightarrow  & J_2(T) & \rightarrow & J_1(F_0) & \rightarrow & Q_1& \rightarrow 0  \\
  &   &  &  \downarrow & & \downarrow & & \downarrow &  \\
  &  &  &  0 & & 0 & &0 &
  \end{array}  \]
and, similarly but after $r\geq 0$ prolongations, the long exact snake-type connecting sequence:  \\
\[  0 \rightarrow R_{r+2} \rightarrow R_{r+1}  \rightarrow h_{r+1} \rightarrow Q_{r+1} \rightarrow Q_r \rightarrow 0  \]
As $Q_1=0$ by counting the dimensions with $dim(R_2)=n+n(n-1)/2=n(n+1)/2$ and $g_3=0$, we get $dim(Q_2) \leq dim(h_2)=n^2(n+1)^2/4 - n^2(n+1)(n+2)/6=n^2(n^2-1)/12$. Hence, we understand that the number of CC ${\cal{D}}_1$ of ${\cal{D}}$ is equal to the number of components of the Riemann tensor {\it if and only if} $R_2$ is formally integrable, that is {\it if and only if} $\omega$ has constant Riemannian curvature, a result first found by L.P. Eisenhart in 1926 ([9]) though in a different setting (See [23] for an explicit modern proof). Such a necessary condition for constructing an exact differential sequence could not have been used by Lanczos because the Spencer "machinery " has only been known after 1970 ([23]). Otherwise, if the metric does not satisfy this condition, CC may exist by using the Petrov classification but have no link with the Riemann tensor. We may therefore define the {\it model vector bundle} $F_1$ with $dim(F_1)=n^2(n^2-1)/12$ in the sense of Lanczos by the short exact sequence:  \\
\[  0  \rightarrow  S_3T^* \otimes T  \rightarrow  S_2T ^*\otimes F_0   \rightarrow  F_1  \rightarrow 0  \]
A result leading to the operator ${\cal{D}}_1:F_0 \stackrel{j_2}{\longrightarrow} J_2(F_0) \rightarrow F_1$ and the:  \\

\noindent
$\bullet$ {\it Second confusion}: Aplying the Spencer operator $\delta$ to the top line of the preceding diagram, a diagonal chase that we have done in many books ([23-26]) or papers ([ 33-36]) we discover that $F_1$ is just the Spencer $\delta$-cohomology $H^2(g_1)$ at ${\wedge}^2T^*\otimes g_1$ along the following short exact sequence:  \\
\[    0 \rightarrow F_1  \rightarrow {\wedge}^2T^*\otimes g_1 \stackrel{\delta}{\longrightarrow}{\wedge}^3T^*\otimes T  \rightarrow 0  \]
because $g_2=0$ and we get the {\it striking formula} where the $+$ signs have been replaced by $-$ signs:  \\
\[ dim(F_1)=n^2(n-1)^2/4 -n^2(n-1)(n-2)/6=n^2(n^2-1)/12\]  
This result, first found in 1978 ([23]) but never acknowledged, clearly exhibit the two well known algebraic properties of the Riemann tensor. We now understand that Lanczos had in mind to linearize the Riemann tensor over the Minkowski metric, {\it exactly like in GR}, in order to construct a Lagrangian as a function of the corresponding linearization $R^k_{l,ij}$ of the Riemann tensor ${\rho}^k_{l,ij}$, transforming the usual variational problem into a variational with a differential constraint described by the Bianchi identities leading to the operatpor ${\cal{D}}_2$. As an equivalent alternative approach, {\it his idea was to consider the curvature as a field by itself} and construct the lagrangian on this field like in EM while adding the Bianchi identities as a differential constraint by using as many Lagrange multiplier as the number of Bianchi identities, a number not known by combinatorics at the time Lanczos was wriring, a result leading to the:  \\

\noindent
$\bullet$ {\it Third confusion}: Lanczos, who also knew continuum mechanics as an engineer, just copied the way used in elasticity (EL) and in electromagnetism (EM), for example introducing a Lagrangian as a function of the deformation $\epsilon=(1/2) \Omega$ while adding a differential constraint described by the vanishing linearized Rieman tensor with therefore as many Lagrange multpliers as the number $n^2(n^2-1)/12$ of components of the Riemann tensor. It is crucial to notice that {\it the same differential sequence is used one step before}, that is with ${\cal{D}}$ and ${\cal{D}}_1$ while he was dealing with ${\cal{D}}_1$ and ${\cal{D}}_2$ previously, that is {\it one step ahead in the sequence}. We have proved recently that such a procedure is in total contradiction with the piezoelectricity and photoelasticity existing between EL and EM (See the picture in [40]).  It thus remains to exhibit the {\it Bianchi} operator exactly as we did for the {\it Riemann} operator, with the same historical comments already provided. However, now we know that $R_1$ is formally integrable (otherwise nothing could be achieved and we should start with a smaller system !), the construction of the linearized Janet-type differential sequence as a strictly exact differential sequence but {\it not} an involutive differential sequence because the system $R_1$ and thus the first order operator ${\cal{D}}$ are formally integrable though {\it not} involutive as $g_1$ is finite type with $g_2=0$ but not involutive. Doing one more prolongation only, we obtain the first order {\it Bianchi} CC from $F_2$ in the following long exact symbol sequence:   \\
\[  0 \rightarrow S_4T^*\otimes T \rightarrow S_3T^*\otimes F_0 \rightarrow T^* \otimes F_1  \rightarrow F_2  \rightarrow 0   \]
or from in the short exact sequence:  \\
\[   0 \rightarrow F_2 \rightarrow {\wedge}^3T^*\otimes g_1 \stackrel{\delta}{\longrightarrow} {\wedge}^4\otimes T  \rightarrow 0  \]
showing that $F_2=H^3(g_1) $, a result still not acknowledged today ([GB1],[GB2], p ...). We have in particular for $n\geq 3$:   \\
\[ \begin{array}{ lcl}
dim(F_2) & = & n^2(n-1)^2(n-2) /12 - n^2(n-1)n-2)(n-3)/24  \\
     & = &  n^2(n+1)(n+2(n+3)/24 +n^3(n^2_1)/12 - n^2(n+1)(n+2)(n+3)/24  \\
& = & n^2(n^2 - 1)(n-2)/24
\end{array}    \]
and thus $dim(F_2)=(4\times 6)-(1\times 4)=(16\times 15 \times 2))24=20$ when $n=4$, a result leading to:   \\

\noindent
$\bullet$ {\it Fourth confusion}: ({\it Double Hodge duality}) For an arbitrry $n$, {\it it is not possible to recognize} that one of the algebraic conditions for the {\it Bianchi} identity comes from the Spencer $\delta$-map and is again an epimorphism as it was before for defining $F_1$, a result obtained by chasing in the commutative diagram obtained by applying $\delta$ to the long exact symbol sequence finishing with $F_2$. It is {\it not evident at all} to discover that the modern description of the model vector bundle $F_2$ is just equivalent to the one provided by Lanczos {\it but only for} $n=4$. For this, using local coordinates, we have the $4$ linear equations with $i=1,2,3,4$:   \\
\[  B^i_{1,234} - B^i_{2,341} + B^i_{3,412} - B^i_{4,123}=0   \]
to be compared with the $4$ equations for the Lanczos tensor $L_{ij,k}=- L_{ji,k}$:   \\
\[  L_{ij,k} + L_{jk,i} + L_{ki,j}=0  \]
Before reading the next lemma, we invite the reader to prove by himself that they are {\it identical}.  \\

\noindent
{\bf LEMMA 3.7}: These two equations are identical only when $n=4$.  \\

\noindent
{\it Proof}: {\it Using Hodge duality a first time}, we may rewrite the first ones in the form:  \\
\[ B^i_{1,1} + B^i_{2,2} + B^i_{3,3}+ B^i_{4,4}=0 , \hspace{5mm} \forall i=1,2,3,4  \]
Lowering the index $i$ by means of the Eucldean metric for simplicity, setting $i=4$ and {\it using again the Hodge duality}, we get with 
$n=4$ only:  \\
\[   B_{44,4}=0 \hspace{5mm}  \Rightarrow \hspace{5mm}  B_{41,1}+B_{42,2} + B_{43,3}=0   \]
Setting now $B_{41,1}= L_{23,1}$ and so on, we get:  \\
 \[    L_{23,1} + L_{31,2} + L_{12,3}=0  \]
 that is e{\it xactly} the Lanczos formula, a result showing that, {\it for} $n=4$ {\it only}, we have $L\simeq B \in {\wedge}^3T^*\otimes g_1$ are both killed by $\delta$.  \\
 \hspace*{12cm}       Q.E.D.   \\

Using adjoint operators and adjoint bundles while setting $ad(E)={\wedge}^nT^*\otimes E^*$ when $E$ is a vector bundle over $X$ and using the Hodge duality, we obtain the short exact sequences with arrows reversed: \\
\[   0 \leftarrow ad(F_2) \leftarrow {\wedge}^{n-3}T^*\otimes g^*_1 \stackrel{\,\,{\delta}^*}{\longleftarrow} {\wedge}^{n-4}T^*\otimes T^* 
\leftarrow  0   \]
as a way to describe the Lagrange multiplier $\lambda \in ad(F_2)$ in arbitrary dimension. \\
We study separately the cases $n=2$ and $n=3$ ([38],[39],[41]):  \\

\noindent
$n=2$:  We have $F_0=S_2T^*, F_1={\wedge}^2T^*\otimes g_1\simeq g_1 \Rightarrow dim(F_1)=1, F_2=0  $ and there is no need at all for any Lagrange multiplier. \\

\noindent
$n=3$: We have $F_0=S_2T^*, dim(F_1)=6, F_2={\wedge}^3T^*\otimes g_1\simeq g_1\Rightarrow dim(F_2)=3$. As we have proved in ([40]) that $Beltrami=ad(Riemann)=Riemann$ is self adjoint up to a change of basis, the corresponding {\it Bianchi} operator is ({\it badly}) identical to the $div$-type {\it Cauchy} operator, a result leading to the:   \\

\noindent
$\bullet$ {\it Fifth confusion}: The $div$-type operator induced ({\it on the right}) by the {\it Bianchi} operator has strictly nothing to do with the {\it Cauchy} operator (namely {\it ad(Killing) on the left}), contrary to what is still believed in GR. In addition, we have the:   \\

\noindent
$\bullet$ {\it Sixth confusion}: We have proved in ([38],[39],[40]) that the usual {\it Cauchy} {\it stress equations} must be replaced by the {\it Cosserat} {\it couple-stress equations} or, equivalently, that the Janet sequence must be replaced by the Spencer sequence in a coherentway with the couplings existing between EL and EM. It is also important to notice that, in the non-linear framework,  there is no analogue of ${\cal{D}}_2$ in the {\it nonlinear Janet sequence} or of $D_3$ in the {\it nonlinear Spencer sequence}, a {\it redhibitory reason} leading to use only ${\cal{D}}$ and 
${\cal{D}}_1$ or $D_1$ and $D_2$ both with their formal adjoints. \\

As a way to conclude this example, we may say that, for any $n\geq 3$, the operator ${\cal{D}}_1= Riemann$ is parametrizing ${\cal{D}}_2=Bianchi$ while $ad({\cal{D}}_2$ is parametrizing the operator $ad({\cal{D}}_1)$. Nevertheless, according to ([37]), there may exist minimal parametrizations of ${\cal{D}}_2$ with a lower number of potentials equal to $dim(T) + dim(F_1) - dim(F_0)$. \\

\noindent
{\bf EXAMPLE 3.8}: ({\it Weyl/Lanczos problem})  \\
Apart from the $6$ confusions that we have pointed out in the preceding Riemann/Lanczos problem, there is a new additional and crucial one coming from the fact that Lanczos and followers tried to use the same $3$-tensor potential $L$ for the system of conformal Killing equations and the coresponding $CKilling$ operator  $\hat{{\cal{D}}}$. More generally, all conformal concepts will be described with a "hat", using the same method as before in order to provide the strictly exact differential sequence:  \\
\[  0 \rightarrow \hat{\Theta} \rightarrow  T  \stackrel{\hat{\cal{D}}}{\longrightarrow} {\hat{F}}_0   \stackrel{{\hat{\cal{D}}}_1}{\longrightarrow} 
{\hat{F}}_1  \stackrel{{\hat{\cal{D}}}_2}{\longrightarrow} {\hat{F}}_2    \]  
where ${\cal{D}}_1$ is the $Weyl$ operator with generating CC ${\cal{D}}_3$. However, it is only in $2016$ (See [37-40] for more details) that we have been able to recover all these operators and confirm with computer algebra the:  \\

\noindent
{\it Seventh confusion}: The orders of the operators involved highly depend on the dimension when $n\geq 3$ as follows:  \\

\noindent 
 $\bullet$ $n=3$: $ \hspace{2cm} 3 \underset{1}{\longrightarrow} 5 \underset{3}{\longrightarrow }5 \underset{1}{ \longrightarrow }3 \rightarrow 0 $ \\
 $\bullet$ $n=4$: $ \hspace{2cm} 4 \underset{1}{\longrightarrow} 9 \underset{2}{\longrightarrow} 10 \underset{2}{\longrightarrow} 9 \underset{1}{\longrightarrow} 4  \rightarrow 0  $ \\                                           
 $\bullet$ $n\geq 5$: $ \hspace{2cm}  5 \underset{1}{\longrightarrow} 14 \underset{2}{\longrightarrow} 35 \underset{1}{\longrightarrow} 35 \underset{2}{\longrightarrow} 14 \underset{1}{\longrightarrow} 5 \rightarrow 0  $\\

This result is based on the following technical lemma (See [23-26] and [36] for details):  \\

\noindent
{\bf LEMMA 3.9}: The symbol ${\hat{g}}_1$ defined by the linear equations:  \\
\[   {\hat{\Omega}}_{ij}\equiv {\omega}_{rj}(x) {\xi}^r_i + {\omega}_{ir}(x){\xi}^r_j - \frac{1}{2}{\omega}_{ij}(x){\xi}^r_r=0  \]
does not depend on any conformal factor, is finite type with ${\hat{g}}_3=0, \forall n\geq 3$ and is surprisingly such that ${\hat{g}}_2$ is $2$-acyclic 
for  $n\geq 4$ or even $3$-acyclic when $n\geq 5$.   \\

Accordigly, we may solve completely the parametrization problem by saying that $ad({\cal{D}}_1)$ is generating the CC of $ad({\cal{D}}_2)$, the number of components of the Lagrange multiplier $\lambda$ being $3$ when $n=3$, $9$ when $n=4$, $35$ when $n=5$ and so on. However, such a result that could be obtained by computer algebra does not provide any information at all on the geometric nature of $\lambda$. The idea will be to use the same diagrams as before but now with ${\hat{g}}_2\simeq T^*$ and ${\hat{g}}_3=0$ in order to obtain in any case ${\hat{F}}_0=\{\bar{\Omega}\mid tr(\bar{\Omega})={\omega}^{ij}{\Omega}_{ij}=0 \}$ and separately: \\

\noindent
$\bullet$ $n=3$: $0\rightarrow {\hat{F}}_1 \rightarrow {\wedge}^2T^*\otimes {\hat{g}}_2 \stackrel{\,\,\delta}{\longrightarrow} 
{\wedge}^3T^*\otimes  {\hat{g}}_1  \rightarrow 0, \hspace{3mm} {\hat{F}}_2={\wedge}^3T^*\otimes {\hat{g}}_2  $.\\

\noindent
$\bullet$ $n=4$: $0 \rightarrow Z^2({\hat{g}}_1) \rightarrow {\wedge}^2T^*\otimes {\hat{g}}_1 \stackrel{\,\,\delta}{\longrightarrow}{\wedge}^3T^*\otimes T \rightarrow 0 , \hspace{1mm} \Rightarrow  \hspace{1mm} 0 \rightarrow T^*\otimes {\hat{g}}_2 \stackrel{\delta}{\longrightarrow} Z^2({\hat{g}}_1) \rightarrow {\hat{F}}_1  \rightarrow 0  $.  \\
\[      0 \rightarrow {\hat{F}}_2 \rightarrow {\wedge}^3T^*\otimes {\hat{g}}_2 \stackrel{\delta}{\longrightarrow} {\wedge}^4T^*\otimes {\hat{g}}_1 \rightarrow 0  \]

\noindent
$\bullet$ $n=5$:  ${\hat{F}}_1=H^2({\hat{g}}_1)  \hspace{4mm}  ,  \hspace{3mm}  {\hat{F}}_2= H^3({\hat{g}}_1)  $ \\

We shall only study the case $n=4$ in order to convince the reader that no classical technique can provide the previous results, also getting in mind, as we saw in the previous sections, that even if an operator is formally integrable, its adjoint may not be formally integrable {\it at all}. For this, using the previous definitions, we first obtain:  \\
\[ dim(T)=4, dim({\hat{F}}_0)= 10-1=9, dim({\hat{F}}_1)= ((6\times 7)-(4\times 4)) - (4\times 4)=26-16=10  \]
Now, if the Bianchi-type operator were first order, applying the Spencer $\delta$-map, we should have ${\hat{F}}_2=H^3({\hat{g}}_1)$ with the long exact symbol sequence:   \\
\[   0 \rightarrow S_4T^*\otimes T \rightarrow S_3T^*\otimes {\hat{F}}_0 \rightarrow T^*\otimes {\hat{F}}_1 \rightarrow H^3({\hat{g}}_1) 
\rightarrow 0   \]
At first, we should get by counting the dimensions $dim(H^3({\hat{g}}_1)=(35\times 4)+(4\times 10) - (20\times 9)= 140 + 40 - 180=0$ but we could also use the formula for cocycle and coboundary:    \\
\[ \begin{array}{rcl}
dim (H^3({\hat{g}}_1))  & = & dim(Z^3({\hat{g}}_1) ) - dim (B^3({\hat{g}}_1))   \\
                                     & =  &(dim ({\wedge}^3T^*\otimes {\hat{g}}_1) - dim({\wedge}^4T^*\otimes T))- 
                                     dim({\wedge}^2T^*\otimes {\hat{g}}_2)  \\
                                     &  = & (4\times 7) - (1\times 4))-(6 \otimes 4)   \\
                                     &  =  &  (28 - 4) - 24=0  
\end{array}  \]
Hence, the generating CC of the Weyl tensor are of order $2$ when $n=4$ and we obtain as before the short exact adjoint sequence:  \\
\[  0 \leftarrow ad({\hat{F}}_2) \leftarrow T^*\otimes {\hat{g}}^*_2  \stackrel{\,\,{\delta}^*}{\longleftarrow}{\hat{g}}^*_1  \leftarrow 0  \]

When $n\geq 5$, using finally the long exact symbol sequence:  \\
\[  0 \rightarrow  S_3T^*  \stackrel{\delta}{\longrightarrow} T^* \otimes S_2T^* \stackrel{\delta}{\longrightarrow}  {\wedge}^2T^* \otimes T^* \stackrel{\delta}{\longrightarrow}  {\wedge}^3T^*  \rightarrow  0   \]
we obtain the commutative and exact {\it fundamental diagram III} playing for the {\it Bianchi} operator the analogue of what the previous {\it fundamental diagram II} was playing for the {\it Riemann} operator:  \\

 \[ \small   \begin{array}{rcccccccll}
 & & & & & & & 0 & &\\
 & & & & & & & \downarrow && \\
  & & & & & 0& & \delta (T^* \otimes S_2 T^*) & & \\
  & & & & & \downarrow & & \downarrow  & & \\
   & & & 0 &\rightarrow & Z^3(g_1) & \rightarrow &H^3({\hat{g}}_1) & \rightarrow 0 & \\
   & & & \downarrow & & \downarrow & &  \downarrow  &  &JANET\\
   & 0 &\rightarrow & {\wedge}^2T^*\otimes {\hat{g}}_2 & \stackrel{\delta}{\rightarrow} & Z^3({\hat{g}}_1) & \rightarrow & H^3({\hat{g}}_1) & \rightarrow 0 & \\
    & & & \downarrow & & \downarrow & & \downarrow &  & \\
 0 \rightarrow & \delta (T^*\otimes S_2T^*)&  \rightarrow  & {\wedge}^2T^*\otimes T^* &\stackrel{\delta}{\rightarrow} & {\wedge}^3T^* & \rightarrow & 0 &  & \\
   & & & \downarrow &  & \downarrow & & & & \\
   & & & 0 & & 0 & & &  &\\
   & & & & & & & & &  \\
   &&&& SPENCER &&&&&
   \end{array}  \]
\vspace*{5mm}  \\  
For $n=5$, we have the following fiber dimensions of the vector bundles involved:

 \[   \begin{array}{rcccccccl}
 & & & & & & & 0 & \\
 & & & & & & & \downarrow & \\
  & & & & & 0& & 40 &  \\
  & & & & & \downarrow & & \downarrow  &  \\
   & & & 0 &\rightarrow & 75 & \rightarrow & 75& \rightarrow 0  \\
   & & & \downarrow & & \downarrow & &  \downarrow  & \\
   & 0 &\rightarrow & 50 & \stackrel{\delta}{\rightarrow} & 85& \rightarrow & 35 & \rightarrow 0  \\
    & & & \downarrow & & \downarrow & & \downarrow &   \\
 0 \rightarrow & 40  & \rightarrow   & 50 &\stackrel{\delta}{\rightarrow} & 10 & \rightarrow & 0 &   \\
   & & & \downarrow &  & \downarrow & & &  \\
   & & & 0 & & 0 & & & 
   \end{array}  \] 

\vspace{4mm}

\noindent
{\bf 3)  CONCLUSION}  \\

\noindent
E. Vessiot discovered the so-called {\it Vessiot structure equations} as early as in 1903 and, only a few years later, E. Cartan discovered the so-called {\it Maurer-Cartan structure equations}. Both are depending on a certain number of constants like the single {\it geometric structure constant} of the constant Riemannian curvature for the first and the many {\it algebraic structure constants} of Lie algebra for the second. However, Cartan and followers never acknowledged the existence of another approach which is therefore still totally ignored today, in particular by physicists. Now, it is well known that the structure constants of a Lie algebra play a fundamental part in the Chevalley-Eilenberg cohomology of Lie algebras and their deformation theory. It was thus a challenge to associate the Vessiot structure constants with other homological properties related to systems of lie equations, namely the extension modules determined by Lie operators. As a striking consequence, such a possibility opens a new way to understand and revisit the various contradictory works done during the last fifty years or so by different groups of researchers, using respectively Cartan, Gr\"obner or Janet bases while looking for a modern interpretation of the work done by C. Lanczos from 1938 to 1962. We hope this paper will open a new domain for applying computer algebra and, in any case, will offer a collection of useful test examples.  \\

\newpage

\noindent
{\bf REFERENCES}\\

\noindent
[1] F. BAMPI, G. CAVIGLIA: Third-order Tensor Potentials for the Riemann and Weyl Tensors, Gen. Relat. and Gravitation, 15 (1983) 375-386.  \\ 
\noindent
[2] E. CARTAN: Les Syst\`{e}mes Diff\'{e}rentiels Ext\'{e}rieurs et Leurs Applications G\'{e}om\'{e}triques, Hermann, Paris,1945.  \\
\noindent 
[3]  P. DOLAN, A. GERBER: Janet-Riquier Theory and the Riemann-Lanczos Problem in 2 and 3 Dimensions, 2002, arXiv:gr-gq/0212055.  \\
\noindent
[4] S.B. EDGAR: Nonexistence of the Lanczos Potential for the Riemann Tensor in Higher Dimension, Gen. Relat. and Gravitation, 26 (1994) 329-332. \\
\noindent
[5] S.B. EDGAR: On Effective Constraints for the Riemann-Lanczos Systems of Equations. J. Math. Phys., 44 (2003) 5375-5385. \\ 
 http://arxiv.org/abs/gr-qc/0302014   \\
\noindent
[6] S.B. EDGAR, A. H\"{O}GLUND: The Lanczos Potential for the Weyl Curvature Tensor: Existence, Wave Equations and Algorithms, Proc. R. Soc. Lond., A, 453 (1997) 835-851.  \\
\noindent
[7] S.B. EDGAR, A. H\"{O}GLUND: The Lanczos potential for Weyl-Candidate Tensors Exists only in Four Dimension, General Relativity and Gravitation, 32, 12 (2000) 2307. \\
http://rspa.royalsocietypublishing.org/content/royprsa/453/1959/835.full.pdf   \\ 
 \noindent
[8] S.B. EDGAR, J.M.M. SENOVILLA: A Local Potential for the Weyl tensor in all dimensions, Classical and Quantum Gravity, 21 (2004) L133.\\
 http://arxiv.org/abs/gr-qc/0408071     \\
\noindent
[9] L.P. EISENHART: Riemannian Geometry, Princeton University Press, 1926.  \\
\noindent
[10] M. JANET: Sur les Syst\`{e}mes aux D\'{e}riv\'{e}es Partielles, Journal de Math., 8(3) (1920) 65-151.  \\
\noindent
[11] E.R. KALMAN, Y.C. YO, K.S. NARENDA: Controllability of Linear Dynamical Systems, Contrib. Diff. Equations, 1, 2 (1963) 189-213.\\
\noindent
[12] M. KASHIWARA: Algebraic Study of Systems of Partial Differential Equations, M\'emoires de la Soci\'et\'e Math\'ematique de France 63, 1995, 
(Transl. from Japanese of his 1970 Master's Thesis).\\
\noindent
[13] E. KUNZ: Introduction to Commutative Algebra and Algebraic Geometry, Birkh\"{a}user, Boston, 1985.  \\
\noindent
[14]  C.LANCZOS: A Remarkable Property of the Riemann-Christoffel Tensor in Four Dimensions, Annals of Math., 39 (1938) 842-850. \\
\noindent
[15] C. LANCZOS: Lagrange Multiplier and Riemannian Spaces, Reviews of Modern Physics, 21 (1949) 497-502.  \\ 
\noindent
[16] C. LANCZOS: The Splitting of the Riemann Tensor, Rev. Mod. Phys. 34, 1962, 379-389.  \\
\noindent
[17] C. LANCZOS: The Variation Principles of Mechanics, Dover, New York, 4th edition, 1949.  \\
\noindent
[18] E. MASSA, E. PAGANI: Is the Rieman Tensor Derivable from a Tensor Potential, Gen.Rel. Grav., 16 (1984) 805-816.  \\ 
\noindent
[19] U. OBERST: Multidimensional Constant Linear Systems, Acta Appl. Math., 20 (1990) 1-175.   \\
\noindent
[20] P. O'DONNELL: A solution of the Weyl-Lanczos equations for the Scwarschild Space-time, for the Scwarzschild Space-Time, General Reativity and Gravitation, 36 (2004) 1415-1422.  \\
\noindent
[21] P. O'DONNELL, H. PYE: A Brief Historical Review of the Important Developments in Lanczos Potential Theory, EJTP, 24 (2010) 327-350.  \\
\noindent
[22] V.P. PALAMODOV: Linear Differential Operators with Constant Coefficients, Grundlehren der Mathematischen Wissenschaften 168, Springer, 1970.\\
\noindent
[23] J.-F. POMMARET: Systems of Partial Differential Equations and Lie Pseudogroups, Gordon and Breach, New York, 1978 
(Russian translation by MIR, Moscow, 1983) \\
\noindent
[24]  J.-F. POMMARET: Differential Galois Theory, Gordon and Breach, New York, 1983.  \\
\noindent
[25]  J.-F. POMMARET: Lie Pseudogroups and Mechanics, Gordon and Breach, New York, 1988.\\
\noindent
[26] J.-F. POMMARET: Partial Differential Equations and Group Theory, New Perspectives for Applications, Mathematics and its Applications 293, Kluwer, 1994.\\
http://dx.doi.org/10.1007/978-94-017-2539-2   \\
\noindent
[27] J.-F. POMMARET: Dualit\'{e} Diff\'{e}rentielle et Applications, C. R. Acad. Sci. Paris, 320, S\'{e}rie I (1995) 1225-1230.  \\
\noindent
[28] J.-F. POMMARET: Partial Differential Control Theory, Kluwer, 2001, 957 pp.\\
\noindent
[29] J.-F. POMMARET: Algebraic Analysis of Control Systems Defined by Partial Differential Equations, in Advanced Topics in Control Systems Theory, Lecture Notes in Control and Information Sciences 311, Chapter 5, Springer, 2005, 155-223.\\
\noindent
[30] J.-F. POMMARET: Gr\"{o}bner Bases in Algebraic Analysis: New perspectives for applications, Radon Series Comp. Appl. Math 2, 1-21, de Gruyter, 2007.\\
\noindent
[31] J.-F. POMMARET: Parametrization of Cosserat Equations, Acta Mechanica, 215 (2010) 43-55.\\
\noindent
[32] J.-F. POMMARET: Spencer Operator and Applications: From Continuum Mechanics to Mathematical Physics, in "Continuum Mechanics-Progress in Fundamentals and Engineering Applications", Dr. Yong Gan (Ed.), ISBN: 978-953-51-0447--6, InTech, 2012, Chapter 1, Available from: \\
http://www.intechopen.com/books/continuum-mechanics-progress-in-fundamentals-\\and-engineering-applications   \\
\noindent
[ 33] J.-F. POMMARET: Relative Parametrization of Linear Multidimensional Systems, Multidim. Syst. Sign. Process. (MSSP), Springer, 26 
(2013) 405-437. \\
http://dx.doi.org/10.1007/s11045-013-0265-0  \\
\noindent
[34] J.-F. POMMARET: The Mathematical Foundations of General Relativity Revisited, Journal of Modern Physics, 4 (2013) 223-239.\\
http://dx.doi.org/10.4236/jmp.2013.48A022   \\
\noindent
[35] J.-F. POMMARET: The Mathematical Foundations of Gauge Theory Revisited, Journal of Modern Physics, 5 (2014) 157-170.  \\
http://dx.doi.org/10.4236/jmp.2014.55026    \\
\noindent
[36] J.-F. POMMARET: From Thermodynamics to Gauge Theory: the Virial Theorem Revisited, in " Gauge Theories and Differential geometry ", NOVA Science Publishers, 2015, Chapter 1, 1-44.  \\
http://arxiv.org/abs/1504.04118  \\
\noindent
[37] J.-F. POMMARET: Airy, Beltrami, Maxwell, Einstein and Lanczos Potentials Revisited, Journal of Modern Physics, 7 (2016) 699-728.  \\
http://dx.doi.org/10.4236/jmp.2016.77068   \\
\noindent
[38] J.-F. POMMARET: Deformation Theory of Algebraic and Geometric Structures, Lambert Academic Publisher, (LAP), Saarbrucken, Germany, 2016.  \\
http://arxiv.org/abs/1207.1964  \\
\noindent
[39] J.-F. POMMARET: Why Gravitational Waves Cannot Exist, Journal of Modern Physics, 8,13 (2017) 20122-2158.  \\
http://dx.doi.org/10.4236/jmp.2017.813130   \\
\noindent
[40] J.-F. POMMARET: From Elasticity to Electromagnetism: Beyond the Mirror, \\
http://arxiv.org/abs/1802.02430  \\
\noindent
[41] J.-F. POMMARET, A. QUADRAT: Localization and parametrization of linear multidimensional control systems, Systems \& Control Letters, 
37 (1999,)247-260.  \\
\noindent
[42] M.D. ROBERTS: The Physical Interpretation of the Lanczos Tensor, Il Nuovo Cimento, B110  (1996) 1165-1176.  \\ 
http://arxiv.org/abs/gr-qc/9904006   \\   
\noindent
[43] J.-P. SCHNEIDERS: An Introduction to D-Modules, Bull. Soc. Roy. Sci. Li\`{e}ge, 63 (1994) 223-295.  \\
\noindent
[44] D.C. SPENCER: Overdetermined Systems of Partial Differential Equations, Bull. Amer. Math. Soc., 75 (1965) 1-114.\\
\noindent
[45] H. TAKENO: On the Spintensor of Lanczos, Tensor, 15 (1964) 103-119.  \\
\noindent
[46] E. VESSIOT: Sur la Th\'{e}orie des Groupes Infinis, Ann. Ec. Normale Sup., 20 (1903) 411-451 (Can be obtained from  http://numdam.org).  \\
\noindent
[47] E. ZERZ: Topics in Multidimensional Linear Systems Theory, Lecture Notes in Control and Information Sciences 256, Springer, 2000.\\

\end{document}